\documentclass[11pt]{amsart}
\usepackage{latexsym,amssymb,amsxtra,amsmath,amsfonts}
\usepackage{verbatim}
\usepackage{color}
\newtheorem{theorem}{Theorem}
\newtheorem{lemma}{Lemma}

\newtheorem{corollary}{Corollary}

\theoremstyle{definition}

\newtheorem{remark}[lemma]{Remark}

\newcommand{\beq}{\begin{equation}}
\newcommand{\eeq}{\end{equation}}
\newcommand{\beqa}{\begin{eqnarray}}
\newcommand{\eeqa}{\end{eqnarray}}
\newcommand{\beaa}{\begin{eqnarray*}}
\newcommand{\ben}{\begin{eqnarray*}}
\newcommand{\eaa}{\end{eqnarray*}}
\newcommand{\een}{\end{eqnarray*}}

\def\F{\mathcal{F}}

\def\L{\mathcal{L}}

\def\ZZ{\mathbb{Z}}

\def\CC{\mathbb{C}}

\def\lieh{{\mathfrak{h}}}

\begin{document} 

\title[2-BKP Grassmanian and D-singularities]
{\textbf{ The 2-component BKP Grassmanian and simple singularities of
    type $D$}}

\author{Jipeng Cheng}
\address{School of Mathematics, Chinese University of Mining and
  technology, P.R. China}
\email{chengjp@cumt.edu.cn}

\author{Todor Milanov}
\address{Kavli IPMU (WPI), UTIAS, The University of Tokyo, Kashiwa, Chiba 277-8583, Japan}
\email{todor.milanov@ipmu.jp}

\begin{abstract}
It was proved in 2010 that the principal Kac--Wakimoto hierarchy
of type $D$ is a reduction of the 2-component BKP hierarchy. On the
other hand, it is known that the total descendant potential of a
singularity of type $D$ is a tau-function of the principal Kac--Wakimoto
hierarchy. We find explicitly the point in the Grassmanian of the
2-component BKP hierarchy (in the sense of Shiota) that corresponds to the total
descendant potential. We also prove that the space of tau-functions of
Gaussian type is parametrized by the base of the miniversal unfolding
of   the simple singularity of type $D$.
\end{abstract}
\maketitle

\section{Introduction}

The total descendant potential in singularity theory is
defined through K. Saito's theory of primitive forms (see \cite{SaK,SaM,He}) and
Givental's higher genus reconstruction \cite{G1}. In the case of
singularities of type $A_{N-1}$ it was proved by Givental in \cite{G3}
that the total
descendant potential is a tau-function of the $N$-KdV hierarchy satisfying
the string equation. Recalling the results of Kac and Schwartz (see
\cite{KS}) such a tau-function is unique and it can be identified with
a point in the big cell of Sato's Grassmanian (see \cite{DJKM, Sa} for
some background)
\ben
\operatorname{Gr}^{(0)} = \{U\subset \CC(\!(z^{-1})\!)\ |\ \pi|_U:U\to
\CC[z] \mbox{ is an isomorphism}\},
\een
where $\CC(\!(z^{-1})\!)$ is the space of formal Laurent series near
$z=\infty$ and $\pi : \CC(\!(z^{-1})\!)\to \CC[z]$ is the linear map
that truncates the terms in the Laurent series that involve negative
powers of $z$. Every point $U$ in Sato's
Grassmanian is determined uniquely by its {\em wave} or {\em Baker
  function}, which by definition is the unique formal function of the type
\ben
\Psi(x,z) = \Big(1+ w_1(x) z^{-1} +w_2(x) z^{-2}+\cdots \Big) e^{xz},\quad
w_i(x)\in \CC[\![x]\!]
\een
such that the Taylor's series expansion in $x$ has coefficients that
span the subspace $U$.  
 Kac and Schwartz proved that the wave function
corresponding to a solution of the $N$-KdV hierarchy satisfying the string
equation can be identified uniquely with a certain solution to the
generalized Airy equation. 

The $N$-KdV hierarchy has a generalization to any simple Lie algebra
of type $ADE$ provided by the so called Kac--Wakimoto hierarchies (see
\cite{KW}). Our main focus will be on the case $D$ and our main goal is to
obtain the analogue of Kac and Schwartz's result. The starting point
of this paper is an observation due to  Liu--Wu--Zhang \cite{LWZ} that
the principal Kac--Wakimoto hierarchy of type $D$ is a reduction of
the 2-component BKP hierarchy. Furthermore, our work relies on the constructions of
Shiota (see \cite{Shi}) and ten Kroode--van de Leur (see \cite{KvL}). 

\subsection{The Grassmanian of the 2-component BKP hierarchy}
Let us recall the construction due to Shiota \cite{Shi} of an inifinte
Grassmanian that plays a key role in the study of the 2-component BKP
hierarchy. Let  
\ben
V:=\CC(\!(z_1^{-1})\!) \oplus \CC(\!(z_2^{-1})\!),
\een
where $z:=(z_1,z_2)$ are formal variables and $( \ ,\ )$ be a symmetric non-degenerate bi-linear pairing defined
by 
\ben
(f(z),g(z)) := \operatorname{Res}_{z_1=0} \frac{dz_1}{z_1} f_1(z_1) g_1(-z_1)
+ \operatorname{Res}_{z_2=0} \frac{dz_2}{z_2} f_2(z_2) g_2(-z_2),
\een
where $f(z)=(f_1(z_1),f_2(z_2)), g(z)=(g_1(z_1),g_2(z_2))\in V$ and the
residues are understood formally as the coefficients in front of
$dz_i/z_i$, $i=1,2$.
The vector space $V$ has a direct sum decomposition $V=U_0\oplus V_0$,
where 
\ben
U_0:=\CC( e_1 +\mathbf{i} e_2) + \CC[z_1] z_1\, e_1 + \CC[z_2]z_2\, e_2,
\een
and 
\ben
V_0:= \CC( e_1 -\mathbf{i} e_2) + \CC[\![z_1^{-1}]\!] z_1^{-1}\, e_1 + \CC[\![z_2^{-1}]\!]z_2^{-1}\, e_2,
\een
where $e_1:=(1,0)$, $e_2:=(0,1)$, and $\mathbf{i}:=\sqrt{-1}$. Let
$\pi:V\to U_0$ be the projection along $V_0$. The big cell of
the Grassmanian of the 2-component BKP hierarchy is by definition the set 
$\operatorname{Gr}_2^{I,(0)} $ of all linear subspaces $U\subset V$ such that
\begin{enumerate}
\item[(1)]
The projection  
$
\pi|_U : U\to U_0
$ 
is an isomorphism.
\item[(2)] U is maximally isotropic.
\end{enumerate}
Recall that a subspace $U$ is called isotropic if $(u_1,u_2)=0$ for
all $u_1,u_2\in U$. It is maximally isotropic if it is not contained
in a larger isotropic subspace, i.e., if $u'\notin U$ then there exist
$u\in U$, such that $(u,u')\neq 0$. Both $U_0$ and $V_0$ are maximal
isotropic subspaces. 

\subsection{The 2-component BKP hierarchy}\label{sec:2-BKP}
Let us recall the correspondence between the points of 
$\operatorname{Gr}_2^{I,(0)}$ and the tau-functions of the 2-component
BKP hierarchy.  Suppose that
\ben
\mathbf{t}^a:=(t^a_m)_{m\in \ZZ_{\rm odd}^+}, \quad a=1,2
\een
are two sequences of formal variables. We denote by $\ZZ_{\rm odd}$
the set of all odd integers and by $\ZZ_{\rm odd}^+$ the set of all
positive odd integers. A formal power series 
\ben
\tau(\mathbf{t}^1,\mathbf{t}^2)\in \CC[\![\mathbf{t}^1,\mathbf{t}^2]\!] 
\een
is said to be a {\em tau-function} of the 2-component BKP if the
following Hirota bilinear equations hold
\ben
\Omega_0(\tau\otimes \tau)=0,
\een
where $\Omega_0$ is the following bi-linear operator acting on
$\CC[\![\mathbf{t}^1,\mathbf{t}^2]\!] ^{\otimes 2}$
\ben
\operatorname{Res}_{z_1=0} \frac{dz_1}{z_1} \Big(
  \Gamma(\mathbf{t}^1,z_1) \otimes \Gamma(\mathbf{t}^1,-z_1) \Big)
  -
\operatorname{Res}_{z_2=0} \frac{dz_2}{z_2} \Big(
  \Gamma(\mathbf{t}^2,z_2) \otimes \Gamma(\mathbf{t}^2,-z_2) \Big),
\een
where 
\ben
\Gamma(\mathbf{t}^a, z_a) := \exp\Big(
\sum_{m\in \ZZ_{\rm odd}^+} t^a_m (z_a)^m\Big)
\exp\Big(-
\sum_{m\in \ZZ_{\rm odd}^+} 2\partial_{t^a_m} \frac{(z_a)^{-m}}{m}\Big)
\een
are {\em vertex operators}. 
Suppose that $\tau$ is a tau-function of the 2-component BKP hierarchy. We will
be interested only in tau-functions such that $\tau(0)\neq 0$, i.e.,
both $\log \tau$ and $\tau^{-1}$ exist in
$\CC[\![\mathbf{t}^1,\mathbf{t}^2]\!]$. The
first two dynamical variables play a special role, so let us denote
them by $x_1=t^1_1$ and $x_2=t^2_1$. We will refer to them as {\em spatial
variables}. Then we define the {\em wave
  function} corresponding to the tau-function $\tau$ to be
\beq\label{wave}
\Psi(x,z) := \Psi^{(1)}(x,z_1) e_1 + \mathbf{i}\, \Psi^{(2)}(x,z_2) e_2,  
\eeq 
where $x=(x_1,x_2)$, $z=(z_1,z_2)$ and the components
\ben
\Psi^{(a)}(x,z_a) := 
\left. 
\frac{\Gamma(\mathbf{t}^a,z_a)\tau(\mathbf{t}^1,\mathbf{t}^2)}{\tau(\mathbf{t}^1,\mathbf{t}^2)}
\right|_{t^b_m=0 (b=1, 2; m>1)},
\quad a=1,2.
\een
The wave function has the following three properties.
\begin{enumerate}
\item[(W1)] The components of the wave function are formal asymptotic series of the type
\ben
\Psi^{(a)}(x,z_a) = \Big(1+\sum_{k=1}^\infty W^{(a)}_k(x) (z_a)^{-k}
\Big) \, e^{x_a z_a},\quad a=1,2,
\een
where $W^{(a)}_k(x)\in \CC[\![x_1,x_2]\!]$. 
\item[(W2)] There exists a formal function $q(x)\in
  \CC[\![x_1,x_2]\!]$ such that 
\ben
(\partial_1\partial_2 +q(x)) \Psi(x,z) = 0,
\een
where $\partial_a=\frac{\partial}{\partial x_a}$ ($a=1,2$). 
\item[(W3)] The wave function is isotropic 
\ben
(\Psi(x',z),\Psi(x'',z)) = 0,
\een
where $x'$ and $x''$ are two copies of the spatial variables
$x=(x_1,x_2)$. 
\end{enumerate}
It turns out that the wave function determines the corresponding
tau-function uniquely up to a constant factor. Moreover, every function of
the form \eqref{wave} satisfying properties (W1)--(W3) is a wave
function, i.e., it corresponds to a tau-function of the 2-component
BKP hierarchy. 

Suppose now that $\Psi(x,z)$ satisfies properties
(W1)--(W3). According to Shiota (see \cite{Shi}, Lemma 9 and its
Corollary) the coefficients of the Taylor's series expansion of
$\Psi(x,z)$ at $x=(0,0)$ span a subspace that belongs to the
Grassmanian $\operatorname{Gr}_2^{I,(0)}$. Moreover, this map
establishes a one-to-one correspondance between the set of wave
functions of the 2-component BKP hierarchy and the points of
$\operatorname{Gr}_2^{I,(0)}$. For proofs of the statements in this
section and for more details we refer to Section 3.1 in \cite{Shi}.

\subsection{Boson-fermion isomorphism}
Let us denote by $Cl(V)$ the  Clifford algebra associated to the vector
space $V$ and the pairing $(\ ,\ )$, i.e., $Cl(V)=T(V)/I$ where 
\ben
T(V):=\CC \oplus \bigoplus_{n\geq 1} V^{\otimes n} 
\een 
is the tensor algebra and $I\subset T(V)$ is the two-sided ideal
generated by  
\ben
v_1\otimes v_2 +v_2\otimes v_1 - (v_1,v_2),\quad v_1,v_2\in V.
\een
Following ten Kroode and van de Leur (see \cite{KvL}) we construct the principal
realization of the basic representation of type $D$ via the spin
representation of $Cl(V)$. Namely, let us introduce the fermionic Fock
space 
\ben
\F := Cl(V)/Cl(V) U_0,
\een
where $U_0\subset V$ is the maximal isotropic subspace introduced
above. The image of $1\in Cl(V)$ in $\F$ will be denoted by
$|0\rangle$ and it will be called the {\em vacuum}. Note that $\F$ is
an irreducible $Cl(V)$-module. We denote by $\phi_a(k)$ the linear
operators induced by multiplication by $e_a (-z_a)^{-k}$. Recalling
the definition of the pairing we get that these operators satisfy the
following relations
\ben
\phi_a(k) \phi_b(\ell)+\phi_b(\ell)\phi_a(k) = (-1)^k\, \delta_{k,-\ell} \delta_{a,b}.
\een
Linear operators satisfying such commutation relations are also known
as {\em neutral free fermions}. 

More explicitly, we can uniquely recover the $Cl(V)$-module structure
on $\F$ by the relations 
\ben
\phi_a(k)\, |0\rangle =0\quad \mbox{for all } k<0, \\
(\phi_1(0)+\mathbf{i}\phi_2(0))\, |0\rangle =0,\\
\phi_1(0)^2 = \phi_2(0)^2 = 1/2 
\een
and the fact that the following set of vectors 
\beq\label{basis-F}
\phi_1(k^1_1)\cdots \phi_1(k^1_r) \phi_2(k^2_1)\cdots \phi_2(k^2_s)\,|0\rangle 
\eeq
where $k^1_1>\cdots >k^1_r >0$ and $k^2_1>\cdots >k^2_s\geq 0$ is a
linear basis of $\F$. 

Let $\F_0$ (resp. $\F_1$) be the subspace of $\F$ spanned by vectors \eqref{basis-F}
such that $r+s$ is even (resp. odd). We would like to equip both
$\F_0$ and $\F_1$ with the structure of an irreducible highets weight
module over a certain Heisenberg algebra. This is a standard
construction. Namely, put 
\ben
J^a_m:=\sum_{j\in \ZZ} (-1)^{j}\, :\phi_a(-j-m)\phi_a(j):,\quad m\in
\ZZ_{\rm odd},\quad a=1,2, 
\een
where the normal ordering is defined by 
\ben
: a b: = ab -\langle a b \rangle,\quad a,b\in V,
\een
where the vacuum expectation $\langle ab \rangle$ is the coefficient
in front of the vacuum $|0\rangle$ when the vector $a b |0\rangle$
is written as a linear combination of the basis vectors \eqref{basis-F}.
The operators satisfy Heisenberg commutation relations
\ben
[J^a_k,J^b_\ell]=2k \delta_{k,-\ell}\delta_{a,b}.
\een
Moreover, the fermions can be expressed in terms of the operators
$J^a_m$ as follows
\ben
\phi_a(z_a):=\sum_{k\in \ZZ} \phi_a(k) (z_a)^k = Q_a 
\exp\Big( \sum_{m\in \ZZ^+_{\rm odd} }J^a_{-m} \frac{(z_a)^m}{m}\Big)
\exp\Big( -\sum_{m\in \ZZ^+_{\rm odd}} J^a_{m} \frac{(z_a)^{-m}}{m}\Big),
\een
where $Q_a:\F\to \F$ are linear operators defined by the following
relations
\ben
Q_a\,|0\rangle & = & \phi_a(0)\, |0\rangle,\quad a=1,2,\\
Q_a \phi_a(k) & = & \phi_a(k) Q_a,\quad a=1,2,\\
Q_a \phi_b(k) & = & -\phi_b(k) Q_a,\quad \mbox{ for } a\neq b.
\een
Finally, there is a unique isomorphism $\F_0\cong
\CC[\mathbf{t}^1,\mathbf{t}^2]$ such that the vacuum $|0\rangle \mapsto
1$ and 
\ben
J^a_{-m}\mapsto m t^a_m,\quad J^a_m\mapsto 2\frac{\partial}{\partial
  t^a_m},\quad 1\leq a\leq 2,\quad m\in \ZZ_{\rm odd}^+.
\een

\subsection{The Kac--Wakimoto hierarchy and the
  2-component BKP hierarchy}
To begin with note that under the boson-fermion isomorphism the
bilinear operator $\Omega_0$ of the Hirota bilinear equations
satisfies the following relations
\ben
\Omega_0=4(Q_1\otimes Q_1)\,\sum_{a=1,2}\sum_{k\in \ZZ} (-1)^k \phi_a(k)\otimes \phi_a(-k) 
,
\een 
where we used that 
\ben
(Q_1\otimes Q_1)^2 = 1/4,\quad (Q_1\otimes Q_1)  (Q_2\otimes Q_2) = -1/4.
\een
Put $h_1:=h:=2N-2$ and $h_2:=2$ ($N\geq 3$). We will be interested in the
so-called $(h_1,h_2)$-reduction of the 2-component BKP. Let us introduce
the bilinear operators $\Omega_m$ ($m\in \ZZ$) acting on 
$\CC[\mathbf{t}^1,\mathbf{t}^2]^{\otimes 2}$ as follows:
\ben
\Omega_m := 4(Q_1\otimes Q_1)\,\sum_{a=1,2}\sum_{k\in \ZZ} (-1)^k \phi_a(-k)\otimes \phi_a(-k-m\,h_a) .
\een
The key observation is that if we set 
\ben
H_{i,m} := \frac{1}{\sqrt{2}} J^1_m,\quad \mbox{for}\quad 1\leq i\leq N-1, \quad m\equiv 2i-1
(\operatorname{mod} h),
\een
and
\ben
H_{n,m(N-1)}:= \sqrt{\frac{N-1}{2}} J^2_m,\quad m\in \ZZ_{\rm odd},
\een
then we get an isomorphism of the Heisenberg algebra spanned by
$J^a_m$ and the principal Heisenberg 
algebra of the affine Lie algebra $\widehat{so}_{2N}$. Therefore, we
can identify $\F_0$ with the principal realization of the basic
representation of $\widehat{so}_{2N}$. In particular, the Casimir of
the Kac--Wakimoto hierarchy acts on $\F_0^{\otimes 2}$ via a certain
bi-linear operator $\Omega_{\rm KW}$.  Our first result can be stated
as follows
\begin{theorem}\label{t1}
The Casimir of the Kac--Wakimoto hierarchy satisfies the following
relation
\ben
-4\Omega_{\rm KW} = \frac{1}{2}\, \Omega_0^2 +\sum_{m=1}^\infty \Omega_{-m}\Omega_m.
\een
\end{theorem}
Note that the action of the operators $\Omega_m$ $(m\in \ZZ)$ extend
to the completion $\CC[\![\mathbf{t}^1,\mathbf{t}^2]\!]$. The above
formula yields the following corollary.
\begin{corollary}\label{c1}
A formal power series $\tau\in
\CC[\![\mathbf{t}^1,\mathbf{t}^2]\!]$ is a solution to the
Kac--Wakimoto hierarchy $\Omega_{\rm KW}(\tau\otimes \tau)=0$ if and
only if one of the following two equivalent conditions are satisfied

a) $\Omega_m(\tau\otimes\tau)=0$ for all $m\geq 0$. 

b) $\tau$ is a tau-function of the 2-component BKP hierarchy and the
corresponding point $U\in \operatorname{Gr}_2^{I,(0)}$ has the
following symmetry
\ben
(z_1^h,z_2^2)\, U\subset U,
\een
where $(z_1^h,z_2^2)$ acts on $V$ by component-wise multiplication
\ben
(z_1^h,z_2^2)\, (f_1(z_1),f_2(z_2)) := (z_1^h f_1(z_1),z_2^2 f_2(z_2)).
\een
\end{corollary}
\begin{remark}
The result in part a) of Corollary \ref{c1} is stated without proof in
\cite{LWZ}. The authors only made a comment that the proof follows from the
fact that the Casimir of the Kac--Wakimoto hierarchy has the form
$\Omega_0^*\Omega_0$. This statement is conceptually true, but
nevertheless a small correction is needed. The goal of Theorem
\ref{t1} is to clarify the remark of Liu--Wu--Zhang. 
\end{remark}
\begin{remark}
Recall that the $N$-KdV hierarchy is a reduction of the KP
hierarchy. It is well known that it can be identified also with
the principal Kac--Wakimoto hierarchy of type $A_{N-1}$. The precise
relation between the Casimirs of the Kac--Wakimoto and the KP
hierarchies  should be given by a formula similar to the one in
Theorem \ref{t1}. Quite surprisingly such a formula seems to be missing in the
literature, or at least we could not find it. 
\end{remark}

\subsection{Virasoro Constraints and the dilaton equation}
According to Frenkel--Givental--Milanov (see \cite{FGM,GM}) the total
descendant potential of the simple singularity of type $D_N$ is a
tau-function of the principal Kac--Wakimoto hierarchy of type
$D_N$. Recalling Corollay \ref{c1} we get that the total descendant
potential is a tau-function 
of the $(h_1,h_2)$-reduction of the 2-component BKP with $h_1=2N-2$ and $h_2=2$. On the other hand,
the total descendant potential is known to satisfy Virasoro
constraints and the {\em dilaton equation} (see \cite{G2}). In the
case at hands the Virasoro 
constraints and the dilaton equation can be stated as follows. 
The Virasoro operators are given by 
\beq\label{vir-L_k}
L_k = -\frac{\mathbf{i}}{2}\, J^1_{1+(1+k)h} + D_k,\quad k\in \ZZ,
\eeq
where
\beq\label{vir-D_k}
D_k=\delta_{k,0}\,\frac{N(h+1)}{24h} + 
\sum_{a=1,2}\frac{1}{4h_a}
\sum_{m\in \ZZ_{\rm odd}}
:J^a_m J^a_{-m+kh_a}:\ .
\eeq
Let us point out that the operator $L_k$ is obtained from $D_k$ via
the so-called {\em dilaton shift} $J^1_{-2N+1}\mapsto J^1_{-2N+1}
-\mathbf{i}\,h$. 

The total descendent potential depends on an extra parameter
$\hbar$. Namely
\beq\label{total-dp}
\mathcal{D}(\hbar,\mathbf{t})=\exp\Big(
\sum_{g=0}^\infty F^{(g)}(\mathbf{t}) \hbar^{g-1}
\Big).
\eeq
In order to identify $\mathcal{D}$ with a tau-function of the
Kac--Wakimoto hierarchy we put $\hbar=1$. 
The dependence on $\hbar$ 
can be recovered via the so called dilaton equation 
\beq\label{dil-eq}
\Big(-\frac{\mathbf{i}h}{h+1} \partial_{t_{1+h}^1}+
\sum_{a=1,2}\sum_{m\in \ZZ^+_{\rm odd}} 
t^a_m \partial_{t^a_m} +
\frac{N}{24} +2\hbar \partial_\hbar\Big)\mathcal{D}=0.
\eeq
Note that the first term is obtained from the second one via the
dilaton shift.  We will say that a formal power series
$\tau(\mathbf{t})$ satisfies the dilaton constraint if the solution to
the linear PDE \eqref{dil-eq} with initial
condition $\mathcal{D}(1,\mathbf{t})=\tau(\mathbf{t})$ is a formal
power series $\mathcal{D}(\hbar,\mathbf{t})$ that
has the form \eqref{total-dp}. 

The Virasoro symmetries can be expressed in terms of
the Grassmanian $\operatorname{Gr}_2^{I,(0)}$ as follows. 
Let us introduce the following differential operators acting
component-wise on $V$
\beq\label{vir-ell}
\ell_k(z) = (\ell^{(1)}_k(z_1),\ell^{(2)}_k(z_2) ),\quad k\in \ZZ,
\eeq
where 
\ben
\ell_k^{(a)}(z_a) := -\mathbf{i} (z_a)^{1+(1+k)h_a}\delta_{1,a}+ \frac{k}{2}
(z_a)^{kh_a} + \frac{1}{h_a} (z_a)^{1+kh_a} \frac{\partial}{\partial
  z_a},\quad a=1,2.
\een
We will prove that if a tau-function of the Kac--Wakimoto hierarchy
satisfies the Virasoro constraints 
$L_k\tau=0$ for all $k\geq -1$ then the corresponding plane
$U\in \operatorname{Gr}_2^{I,(0)}$ satisfies 
\ben
\ell_k(z) U \subset U,\quad k\geq -1.
\een
Note that the constraint for $k=-1$, also known as the
{\em string equation}, and the symmetry
$(z_1^h,z_2^2)U\subset U$ imply the rest of the Virasoro constraints. 
The main result of this paper can be stated as follows. 
\begin{theorem}\label{t2}
If a tau-function of the Kac--Wakimoto hierarchy satisfies the
Virasoro and the dilaton constraints, then the
corresponding wave function $\Psi(x,z)$ satisfies the following system
of PDEs 
\begin{align}
\notag
(\partial_1^h + \partial_2^2 - \mathbf{i}\, x_1) \Psi(x,z) & = 
(z_1^h,z_2^2) \Psi(x,z) ,\\
\notag
\partial_1 \Psi(x,z) & =  \mathbf{i} \ell_{-1}(z) \Psi(x,z),\\
\notag
\partial_1\partial_2 \Psi(x,z) & =  \frac{\mathbf{i}}{2} x_2 \Psi(x,z),
\end{align}
where $\partial_1:=\partial/\partial x_1$,
$\partial_2:=\partial/\partial x_2$.
\end{theorem}

\begin{corollary}\label{c2}
The total descendant potential is the unique tau-function of the
Kac--Wakimoto hierarchy satisfying the string and the dilaton
constraints. 
\end{corollary}
\begin{remark}
In the case of singularities of type $A$ the Virasoro constraints
uniquely determine the tau-function. We do not know whether the
additional dilaton constraint imposed in Corollary \ref{c2} is
necessary. 
\end{remark}

\subsection{Tau-functions of Gaussian type}
The problem of classifying tau-functions whose logorithm is a
quadratic form in the dynamical variables was proposed by Givental in
the settings of the $N$-KdV hierarchy (see \cite{G3}). We solve this
problem in the case of the principal Kac--Wakimoto hierarchy of type
$D$. The answer is quite suggestive. Namely, one can speculate that
there is a general theory of Hirota bilinear equations in which the
tau-functions of Gaussian type provide an embedding of the theory of
semi-simple Frobenius manifolds into the theory of integrable
systems. This expectation is in some sense compatible with the general framework
developped by Dubrovin and Zhang in \cite{DZ}. 

Suppose that 
\beq\label{tau-gaussian}
\tau(\mathbf{t}^1,\mathbf{t}^2) = \exp\Big(\frac{1}{2}
\sum_{a,b=1}^2\sum_{k,\ell\in \ZZ^+_{\rm odd}} W^{ab}_{k\ell} t^a_kt^b_\ell\Big).
\eeq
Under what conditions $\tau$ is a tau-function of the Kac--Wakimoto
hierarchy? The answer is given in terms of the following system of equations
\ben
x^{2(N-1)} +\sum_{i=1}^{N-1} t_i x^{2(N-1-i)} + y^2 & = & \lambda,\\
xy+t_N & = &  0,
\een
where $t=(t_1,\dots,t_N)\in \CC^N$ are fixed parameters. We will prove
that the set of all coefficients $W^{ab}_{k\ell}$ is uniquely
determined from the subset of coefficients with $k=1$. The
coefficients $W^{ab}_{1\ell}$ are determined uniquely from the above
system of algebraic equations as follows. Put $\lambda=z_1^h$. Then
the system admits a formal solution of the form
\ben
x=z_1-2\sum_{\ell} W^{11}_{1\ell}(t) z_1^{-\ell}/\ell,\quad
y=-2\sum_{\ell} W^{21}_{1\ell}(t) z_1^{-\ell}/\ell,
\een
where both sums are over all $\ell \in \ZZ_{\rm odd}^+$ and 
$W^{11}_{1\ell}(t),  W^{21}_{1\ell}(t)\in \CC[t]$. Similarly, put
$\lambda=z_2^2$, then the system admits a formal solution of the form
\ben
x=-2\sum_{\ell} W^{12}_{1\ell}(t) z_2^{-\ell}/\ell,\quad
y=z_2-2\sum_{\ell} W^{22}_{1\ell}(t) z_2^{-\ell}/\ell.
\een
\begin{theorem}\label{t3}
The map $t\mapsto W^{ab}_{k\ell}(t)$ defined above establishes a
one-to-one correspondence between $\CC^N$ and the set of tau-functions
of Gaussian type.
\end{theorem}

{\bf Acknowledgements.}
The work of T.M. is partially supported by JSPS Grant-In-Aid (Kiban C)
17K05193 and by the World Premier International Research Center
Initiative (WPI Initiative),  MEXT, Japan. 

We would like to acknowledge also the fact that the formula in Theorem \ref{t1}
was derived with techniques and ideas developed in the joint work of
T.M. with B. Bakalov in \cite{BM}.  

\section{The principal Kac--Wakimoto hierarchy of type $D_N$}
The main goal of this section is to prove Theorem \ref{t1}. We will
use the language of twisted representations of lattice vertex
algebras, because it seems to be the most appropriate one. We have
verified the formula in Theorem \ref{t1} also directly but the 
computation is very long. For some quick
introduction to twisted vertex algebra modules with an extensive list of references for
more systematic study we  refer to \cite{BM}. 
\subsection{The Euclidean lattice vertex algebra}
Let $\lieh:=\CC^N$ and $\{v_i\}_{i=1}^N\subset \lieh$ be the
standard basis. Put
\ben
\ZZ^N:=\bigoplus_{i=1}^N \ZZ\, v_i,\quad (v_i|v_j)=\delta_{i,j}
\een
for the standard Euclidean lattice. Let 
\ben
\epsilon:\ZZ^N\times \ZZ^N\to \{\pm 1\},
\een
be the bi-multiplicative function defined by 
\ben
\epsilon(v_i,v_j) = 
\begin{cases}
-1 & \mbox{ if $i\leq j$ }, \\
1 & \mbox{ if $i>j$}.
\end{cases}
\een
Recall the {\em twisted group} algebra $\CC_\epsilon[\ZZ^N]$ with
basis $e^a$ ($a\in \ZZ^N$) and multiplication 
\ben
e^a\, e^b:= \epsilon(a,b) e^{a+b}. 
\een 
As a vector space the lattice vertex algebra is defined by 
\ben
V_{\ZZ^N}:= \operatorname{Sym}(\lieh[s^{-1}]s^{-1})\otimes \CC_\epsilon[\ZZ^N],
\een
where $\operatorname{Sym}(A)$ is the symmetric algebra of a vector
space $A$. 
Slightly abusing the notation we will write $e^a$ for $1\otimes e^a$ and
$a$ for $(as^{-1})\otimes e^0$. The vector $\mathbf{1}:=1\otimes e^0$
is called the {\em vacuum}. 

The structure of a vertex algebra is given by a bi-linear map 
\ben
Y(\cdot, z): V_{\ZZ^N}\otimes V_{\ZZ^N}\to V_{\ZZ^N}(\!( z)\!),
\een
which is also known as the {\em state-field correspondence}.  Let us
recall the construction of $Y$. 
Let $\widehat{\lieh}=\lieh[s,s^{-1}]\oplus \CC\, K$ be the Heisenberg
Lie algebra whose commutator is defined by 
\ben
[as^m,bs^n]=m\delta_{m,-n} (a|b) \, K.
\een
There is a unique way to turn $V_{\ZZ^N}$ into a
$\widehat{\lieh}$-module such that $as^m$ for $m<0$ acts as
multipliction by $(as^m)\otimes 1$, the central element $K$ acts by 1,
and for $m\geq 0$ 
\ben
as^m (1\otimes e^b) = \delta_{m,0} (a|b) \, 1\otimes e^b,\quad a\in
\lieh,\quad b\in \ZZ^N.
\een
The linear operator representing $as^n\in \widehat{\lieh}$ is denoted
by $a_{(n)}$. 
Put
\ben
Y(a,z) := \sum_{n\in \ZZ} a_{(n)} z^{-n-1},\quad a\in \lieh
\een
and 
\ben
Y(e^b,z):= e^b z^{b_{(0)}} 
e^{\sum_{j>0} b_{(-j)} \frac{z^j}{j}} 
e^{\sum_{j<0} b_{(-j)} \frac{z^j}{j}} ,\quad b\in \ZZ^N.
\een
The definition of the state-field correspondence can be extended
uniquely so that the following formula holds
\beq\label{va:ope-voa}
Y(a_{(n)} b,z) =\frac{1}{k!} \left.\partial_w^k \Big( (w-z)^{n+1+k} Y(a,w)Y(b,z)\Big)\right|_{w=z}
\eeq
for every $a,b\in V_{\ZZ^N}$, where we choose $k\gg 0$ so big that 
\beq\label{locality}
(w-z)^{n+1+k} [Y(a,w),Y(b,z)]=0,
\eeq
and denote by $a_{(n)}$ the Fourier modes $Y(a,z)=:\sum_{n}
a_{(n)}z^{-n-1}$. This 
formula defines $Y(a_{(n)}b,z)$ assuming that we already
know $Y(a,z)$ and $Y(b,z)$. Clearly this is a recursive procedure that 
defines the state-field correspondence in terms of $Y(a,z)$ ($a\in
\lieh$) and $Y(e^a,z)$ ($a\in \ZZ^N$). For more details we refer
to \cite{BM}, Proposition 3.2. 

The main property of the above definition is the so-called Borcherd's
identity for the modes
\begin{equation}\label{borcherd-id}
\begin{split}
\sum_{j=0}^\infty (-1)^j & \binom{n}{j} 
\Bigl( 
a_{(m+n-j)}(b_{(k+j)}c)
- (-1)^n \, b_{(k+n-j)}(a_{(m+j)}c)
\Bigr)
\\ 
&=
\sum_{j=0}^\infty \binom{m}{j} (a_{(n+j)}b)_{(k+m-j)}c \,,
\end{split}
\end{equation}
where $a,b,c \in V_{\ZZ^N}$. Observe that the above sums are finite, because
$a_{(n)}b = 0$ for sufficiently large $n$. We
will need also the following commutator formula 
\ben
[Y(a,w),Y(b,z)]=\sum_{n=0}^\infty \frac{1}{n!} Y(a_{(n)}
b,w)\partial_w^n \delta(z,w),
\een
where $\delta(z,w)=\sum_{m\in \ZZ} z^m w^{-m-1}$ is the formal
delta-function. The above formula is equivalent to the Borcherd's
identity \eqref{borcherd-id} with $n=0$. For more details on lattice
vertex algebras and for the proofs of the statements in  this section
we refer to \cite{K2}. 

\subsection{The Frenkel-Kac construction}

Let us recall the Frenkel--Kac construction (see
\cite{FK}) of the affine Kac--Moody Lie algebra of type $D_N$. 
Recall that the root system of type $D_N$ consists of the following vectors
\ben
\Delta=\{ \pm (v_i\pm v_j)\ |\ 1\leq i\neq j\leq N\}.
\een
The lattice vertex algebra $V_{\ZZ^N}$ has a {\em conformal vector}
\beq\label{conf-vector}
\nu:=\frac{1}{2} \sum_{i=1}^N v_{i(-1)} v_{i(-1)} \mathbf{1}.
\eeq
The affine Kac--Moody Lie algebra of type
$D_N$ can be constructed by the following simple formula
\ben
\widehat{so}_{2N}(\CC)\cong \bigoplus_{\alpha\in \Delta} \bigoplus_{n\in \ZZ} \CC e^{\alpha}_{(n)}
\oplus \bigoplus_{i=1}^N \bigoplus_{n\in \ZZ} \CC v_{i(n)} \oplus \CC
\nu_{(1)} \oplus \CC \, \mathbf{1}_{(-1)}\ ,
\een
where the RHS is equipped with a Lie bracket given by the commutator: 
\ben
[x_{(m)},y_{(n)}]:=x_{(m)}y_{(n)}-y_{(n)}x_{(m)},\quad x,y \in
V_{\ZZ^N},\quad m,n\in \ZZ.
\een
We leave it to the reader as an exercise to use the Borcherd's
identity \eqref{borcherd-id} with $n=0$ to check that the above
formula indeed defines a Lie algebra isomorphic to
$\widehat{so}_{2N}(\CC)$. 

\subsection{Casimirs of $\widehat{so}_{2N}(\CC)$}
The tensor product $V_\ZZ\otimes V_\ZZ$ also has the stucture of a vertex
operator algebra with state-field correspondence defined by 
\ben
Y(a\otimes b, z) := Y(a,z)\otimes Y(b,z).
\een
Let us define the following two vectors in $V_{\ZZ^N}^{\otimes 2}$:
\ben
\omega_{\rm KW} = -\sum_{\alpha\in \Delta} e^{\alpha}\otimes
e^{-\alpha} + \sum_{i=1}^N v_i\otimes v_i - \nu\otimes 1 - 1\otimes \nu
\een
and 
\ben
\omega_{\rm BKP} = \sum_{i=1}^N (e^{v_i}\otimes e^{-v_i}+e^{-v_i}\otimes e^{v_i}).
\een
\begin{lemma}\label{le:casimir}
The following relation holds
\ben
\omega_{\rm KW} = -\frac{1}{2}(\omega_{\rm BKP} )_{(-1)} \omega_{\rm BKP}.
\een
\end{lemma}
\proof
By definition
\ben
(\omega_{\rm BKP} )_{(-1)} \omega_{\rm BKP} = 
\operatorname{Res}_{z=0} \frac{dz}{z}\Big(
Y(\omega_{\rm BKP},z) \omega_{\rm BKP}\Big).
\een
For $\lambda,\mu\in \{+1,-1\}$ and $1\leq i,j\leq N$ we have
\ben
&&
\frac{dz}{z}\Big( (Y(e^{\lambda v_i},z)\otimes Y(e^{-\lambda v_i},z) ) e^{\mu
  v_j}\otimes e^{-\mu v_j} \Big)= 
\frac{dz}{z} \, z^{2\lambda\mu \delta_{i,j}}\times 
\\
&&
(e^{\sum_{n>0} \lambda v_{i(-n)}\, \frac{z^n}{n}} 
e^{\lambda v_i + \mu v_j})\otimes 
(e^{-\sum_{n>0} \lambda v_{i(-n)}\, \frac{z^n}{n}} 
e^{-\lambda v_i - \mu v_j}).
\een
The residue of the above 1-form is not $0$ only if $i\neq j$ or $i=j$
and $\mu=-\lambda$. In the former case the residue is $e^\alpha\otimes
e^{-\alpha}$ with $\alpha= \lambda v_i + \mu v_j\in \Delta$, while in
the latter it is
\ben
 \frac{\lambda}{2} (v_{i(-2)}\otimes 1 - 1\otimes v_{i(-2)} ) + \frac{1}{2}( v_{i(-1)}^2\otimes 1 +1\otimes
 v_{i(-1)}^2) -v_{i(-1)}\otimes v_{i(-1)}. 
\een
Summing over all $i,j\in \{1,2,\dots,N\}$ and over all $\lambda,\mu
\in \{+1,-1\}$ we get 
\ben
\operatorname{Res}_{z=0} \frac{dz}{z}\Big(
Y(\omega_{\rm BKP},z) \omega_{\rm BKP}\Big) = -2\omega_{\rm KW}.
\een

\subsection{The Coxeter transformation}\label{sec:Coxeter}
Let $\alpha_i=v_i-v_{i+1}$ ($1\leq i\leq N-1$) and
$\alpha_N=v_{N-1}+v_N$ be a set of simple roots. Let us fix a Coxeter
transformation 
\ben
\sigma:=r_{\alpha_1}\cdots r_{\alpha_N},
\een
where $r_{\alpha_i}(x)=x-(\alpha_i|x)\alpha_i$ are the simple
reflections. The action of $\sigma$ on the standard basis is
represented by the diagram
\ben
v_1\mapsto v_2\mapsto \cdots \mapsto v_{N-1}\mapsto -v_1,\quad
v_N\mapsto -v_N.
\een
Let us choose an eigenbasis $\{H_i\}_{i=1}^N$ of $\sigma$ such that 
\ben
\sigma(H_i)=\eta^{m_i}H_i,\quad (H_i|H_j)=h \delta_{i,j^*},
\een
where $\eta:=e^{2\pi \mathbf{i}/h}$ (recall that $h=2N-2$ is the
Coxeter number), $m_i:=2i-1$ ($1\leq i\leq N-1$) and 
  $m_N=N-1$ are the so-called {\em exponents}, and ${}^*$ is an
  involution defined by
\ben
j^*:=
\begin{cases}
n-j & \mbox{ if $1\leq j\leq N-1$}, \\
j & \mbox{if $j=N$}.
\end{cases}
\een
To be more specific, let us define $H_i$ as the solutions to the following
linear system of equations (with unknowns $H_i$):
\begin{align}\notag
v_i & = \frac{\sqrt{2}}{h}
\Big( \eta^{m_1 i}H_1 +\cdots +\eta^{m_{N-1}i}H_{N-1}\Big),\quad 1\leq
      i\leq N-1, \\
\notag
v_N & = \frac{1}{\sqrt{h}} \, H_N.
\end{align}
Following Bakalov--Kac (see \cite{BK}, Proposition 4.1) we extend 
the Coxeter transformation $\sigma$ to a vertex
algebra automorphism of $V_{\ZZ^N}$. The action of $\sigma$ on 
$\operatorname{Sym}(\lieh[s^{-1}]s^{-1})$ is induced from the action
of $\sigma$ on $\lieh$. While the definition of the action of $\sigma$
on the twisted group algebra involves the choice of a function  
\ben
\zeta: \ZZ^N \to \{+1,-1\},
\een
such that
\ben
\epsilon(\sigma(a),\sigma(b)) \epsilon(a,b)^{-1} = \zeta(a+b)
\zeta(a)^{-1} \zeta(b)^{-1},\quad a,b \in \ZZ^N. 
\een
The definition 
\ben
\sigma(e^a) = \zeta(a) e^{\sigma(a)},\quad a\in \ZZ^N
\een
extends uniquely to an automorphism of the twisted group algebra
$\CC_\epsilon[\ZZ^N]$. The linear map 
\ben
\sigma: V_{\ZZ^N}\to V_{\ZZ^N},\quad \sigma(x\otimes e^a):=
\sigma(x)\otimes \sigma(e^a),
\een
is an automorphism of vertex algebras, i.e.,
$\sigma(x_{(n)}y)=\sigma(x)_{(n)}\sigma(y)$ for all $x,y\in V_{\ZZ^N}$
and for all $n\in \ZZ$.

Let us define 
\ben
\pi_n(v):=\frac{1}{h}\sum_{j=1}^h \eta^{jn}\sigma^j(v),\quad v\in V_{\ZZ^N}.
\een
Note that $\pi_n$ is the projection of $v$ onto the
eigensubspace of $\sigma$ corresponding to eigenvalue
$e^{-2\pi\mathbf{i} n/h}$. After a small modification of the
Frenkel--Kac construction we define the following Lie algebra:
\ben
\widehat{so}_{2N}(\CC,\sigma):=
\bigoplus_{\alpha\in \Delta/\sigma} \bigoplus_{n\in \ZZ} \CC (e^{\alpha})^\sigma_{(n)}
\oplus \bigoplus_{i\in\{1, N\}} \bigoplus_{n\in \ZZ} \CC (v_i)^\sigma_{(n)} \oplus \CC
\nu_{(1)} \oplus \CC \, \mathbf{1}_{(-1)}\ ,
\een
where for $x\in V_{\ZZ^N}$ we put
$x^\sigma_{(n)}:=(\pi_n(x))_{(n)}$ and $\Delta/\sigma$ is the set of
orbits of $\sigma$ in $\Delta$. Note that if $x$ and $y$ belong to the
same orbit of $\sigma$, i.e., $x=\sigma^\ell (y)$ for some $\ell\in \ZZ$,
then the complex lines $\CC\, x^\sigma_{(n)}$ and
$\CC\,y^\sigma_{(n)}$ coincide.
Using the Borcherd's identities one can check that
$\widehat{so}_{2N}(\CC,\sigma)$ is isomorphic to the {\em twisted
  affine Lie algebra} of type $D_N$ corresponding to the automorphism
$\sigma$. We refer to \cite{K1}, Section 8.2 for some background on twisted
affine Lie algebras. 

\subsection{Twisted representations}

Suppose that $M$ is a vector space and that
\ben
Y^M(\cdot,\lambda):V_{\ZZ^N}\otimes M\to M(\!(\lambda^{1/h})\!)
\een
is a linear map that defines a $\sigma$-twisted representation of
$V_{\ZZ^N}$ on $M$ (see \cite{BK}, Definition 3.1).  If $a\in
V_{\ZZ^N}$ then we denote by $a^M_{(m)}$ $(m\in \frac{1}{h} \ZZ)$ the
modes of the twisted field
\ben
Y^M(a,\lambda) = \sum_{m\in \frac{1}{h} \ZZ^N} a^M_{(m)}\, \lambda^{-m-1}.
\een
Let us recall the following three properties (see \cite{BM}) of a
twisted representation: 
\begin{enumerate}
\item[(i)]$\sigma$-invariance: 
$Y^M(\sigma(v),\lambda)$ coincides with
the analytic continuation in counter clock-wise direction around
$\lambda=0$ of  $Y^M(v,\lambda)$.
\item[(ii)] Locality: for every $a,b\in V_{\ZZ^N}$ there exists
$n_{ab}\geq 0$ such that 
\ben
(\lambda_1-\lambda_2)^{n_{ab}}[Y^M(a,\lambda_1),Y^M(b,\lambda_2)]=0.
\een
\item[(iii)] Product formula: \eqref{va:ope-voa} remains true if we replace $Y$
with $Y^M$, i.e., 
\beq\label{ope-voa}
Y^M(a_{(n)} b,\lambda) =
\frac{1}{k!} \left.\partial_\mu^k \Big( 
(\mu-\lambda)^{n+1+k} Y^M(a,\mu)Y^M(b,\lambda)\Big)\right|_{\mu=\lambda},
\eeq
for all $a,b\in V_{\ZZ^N}$ and $n\in \ZZ$. 
\end{enumerate}
In fact the above properties and the vacuum axiom
$Y_\sigma(\mathbf{1},\lambda)=1$ can be used as a
definition of a twisted representation. In particular, the
locality and the product formula imply the Borcherd's identity for the
twisted modes
\begin{equation}\label{tw_borcherd-id}
\begin{split}
\sum_{j=0}^\infty (-1)^j & \binom{n}{j} 
\Bigl( 
a^M_{(m+n-j)}(b^M_{(k+j)}c)
- (-1)^n \, b^M_{(k+n-j)}(a^M_{(m+j)}c)
\Bigr)
\\ 
&=
\sum_{j=0}^\infty \binom{m}{j} (a_{(n+j)}b)^M_{(k+m-j)}c \,,
\end{split}
\end{equation}
for all $a\in V_{\ZZ^N}$, s.t., $\sigma(a)=e^{-2\pi\mathbf{i}m}a$,
$b\in V_{\ZZ^N}$, $c\in M$, $m,k\in \frac{1}{h} \ZZ$, and
$n\in \ZZ$.

Using the two Borcherd's identities \eqref{borcherd-id} and 
\eqref{tw_borcherd-id} it is straightforward to check that the
maps 
\ben
(e^\alpha)^\sigma_{(n)} \mapsto (e^\alpha)^M_{(n/h)},\quad (\alpha\in
\Delta,\quad n\in \ZZ),\quad 
\nu_{(1)}  \mapsto h \nu^M_{(1)}, \quad
\mathbf{1}_{(-1)} \mapsto h^{-1} {\rm Id}_M
\een
define a representation of $\widehat{so}_{2N}(\CC,\sigma)$ on $M$. Let us
observe also that the standard Euclidean pairing on $\lieh$ can be
extended uniquely to an invariant bi-linear form on
$\widehat{so}_{2N}(\CC)$, such that 
\ben
(e^\alpha_{(m)}|e^\beta_{(n)})=-\delta_{\alpha,-\beta}\delta_{m,-n},\quad
(a_{(m)}|b_{(n)}) = (a|b)\delta_{m,-n},\quad (\nu_{(1)}|\mathbf{1}_{(-1)})=-1,
\een
where $a,b\in \lieh$, $\alpha,\beta\in \Delta$, $m,n\in \ZZ$, and all
other pairings vanish. The twisted affine Lie algebra
$\widehat{so}_{2N}(\CC,\sigma)$ is a Lie subsalgebra of
$\widehat{so}_{2N}(\CC)$. The induced bi-linear form is still
non-degenerate and invariant. Therefore we can construct a Casimir for
$\widehat{so}_{2N}(\CC,\sigma)$ via
$\sum_{a} x_a \otimes x^a$, where $\{x_a\}$ and $\{x^a\}$ is a pair
of dual bases of $\widehat{so}_{2N}(\CC,\sigma)$.  Let us fix a basis
of the following type
\ben
(e^{\alpha_i})^\sigma_{(n)},\quad 
(v_1)^{\sigma}_{(m)},\quad
(v_N)^\sigma_{(m(N-1))},\quad \nu_{(1)},\quad \mathbf{1}_{(-1)},
\een
where $n\in \ZZ $, $m\in \ZZ_{\rm odd} $, and the roots $\alpha_i\in
\Delta$ ($1\leq i\leq N)$ are such that 
the corresponding orbits of the Coxeter transformation are pairwise
disjoint.  Using the formulas
\begin{align}
\notag
((e^{\alpha_i})^\sigma_{(m)}|(e^{-\alpha_j})^\sigma_{(n)}) & = -\delta_{m,-n}\delta_{i,j}/h\\
((v_i)^\sigma_{(m)} | (v_i)^\sigma_{(n)}) & =
\begin{cases}
\frac{1}{N-1}\, \delta_{m,-n} & \mbox{ if $1\leq i\leq N-1$}, \\
\notag
\delta_{m,-n} & \mbox{ if $i=N$},
\end{cases}
\end{align}
it is straightforward to construct a dual basis:
\ben
-h (e^{-\alpha_i})^\sigma_{(-n)}\quad 
(N-1)(v_1)^{\sigma}_{(-m)}\quad
(v_N)^\sigma_{(-m(N-1))},\quad
-\mathbf{1}_{(-1)}, \quad
-\nu_{(1)}.
\een
The Casimir takes the form
\ben
&&
-\sum_{\alpha\in \Delta}\sum_{n\in \ZZ} (e^\alpha)^\sigma_{(n)}\otimes
(e^{-\alpha})^\sigma_{(-n)} 
-\nu_{(1)}\otimes \mathbf{1}_{(-1)} - \mathbf{1}_{(-1)} \otimes
\nu_{(1)} + \\
&&
+\sum_{i=1}^N \sum_{n\in \ZZ} (v_i)^\sigma_{(n)}\otimes
(v_i)^\sigma_{(-n)}\, ,
\een
where we used that the expressions $(e^\alpha)^\sigma_{(n)}\otimes
(e^{-\alpha})^\sigma_{(-n)}$ and $(v_i)^\sigma_{(n)}\otimes
(v_i)^\sigma_{(-n)}$ are invariant under the Coxeter transformations
$\alpha\mapsto \sigma \alpha$ and $v_i\mapsto \sigma v_i$, respectively. 
Note that the action of the Casimir on $M^{\otimes 2}$ is given by
$\operatorname{Res}_{\lambda=0} 
\left(Y^M(\omega_{\rm  KW},\lambda)\lambda d\lambda\right).$ 

\subsection{Twisted modules over the Euclidean lattice vertex algebra}
The $\sigma$-twisted modules over a lattice vertex algebra are classified by
Bakalov and Kac in \cite{BK}. Let us recall their result in the case
of $V_{\ZZ^N}$ when $\sigma$ is the Coxeter transformation. 

The $\sigma$-twisted module structure is determined
uniquely by the so called {\em Heisenberg pair}
$(\widehat{\lieh}_\sigma, G_\sigma)$, which is defined as follows.
The first member of the pair is the $\sigma$-twisted Heisenberg algebra 
\ben
\widehat{\lieh}_\sigma =\CC\, K \oplus
\bigoplus_{i=1}^N\bigoplus_{n\in \ZZ} \CC\, H_{i} s^{-n-m_i/h}  ,
\een
where the notation is the same as in Section \ref{sec:Coxeter}. The
Lie bracket is given by
\ben
[h's^m,h''s^n]:= m\delta_{m,-n}(h'|h'')\, K,\quad h',h''\in
\lieh,\quad m,n\in \frac{1}{h}\ZZ.
\een
Let $G=\CC^*\times \ZZ^N$ be the set whose elements will be written as
$cU_\alpha$, $c\in \CC^*$, $\alpha\in \ZZ^N$. The following
multiplication turns $G$ into a group
\ben
U_\alpha U_\beta:= \epsilon(\alpha,\beta) B(\alpha,\beta)^{-1} U_{\alpha+\beta},
\een
where 
\ben
B(\alpha,\beta) := h^{-(\alpha|\beta)} \prod_{k=1}^{h-1} (1-\eta^k)^{(\sigma^k(\alpha)|\beta)}. 
\een
Put
\ben
N_\sigma:= \{ \zeta(\alpha)^{-1} U_{\sigma\alpha}^{-1} U_\alpha \,
(-1)^{|\alpha|^2}\ |\ \alpha \in \ZZ^N\},
\een
where $|\alpha|^2:=(\alpha|\alpha)$.
Using the commutator formula 
\ben
U_\alpha U_\beta U_{\alpha}^{-1} U_\beta^{-1}=
\exp\ 2\pi\mathbf{i}\Big(\frac{1}{2} |\alpha|^2 |\beta|^2 + ((1-\sigma)^{-1}\alpha|\beta)\Big)
\een
it is easy to check that $N_\sigma$ is a subroup of the center
$Z(G)$. 
The second member of the Heisenberg pair is the quotient group
$G_\sigma:=G/N_\sigma$. Let us denote by $\overline{U}_\alpha\in
G_\sigma$ the image of $U_\alpha$ under the quotient map. We will be
interested in representations of $G_\sigma$ that are $\CC^*$-invariant, i.e., the element
$c\overline{U}_0\in G_\sigma$ ($c\in \CC^*$) acts by multiplication by
the scalar $c$.
\begin{lemma}
a) The following relations hold
\ben
\overline{U}_{v_1} ^2 = (-1)^{N} \zeta(v_1)\cdots \zeta(v_{N-1})\,
\frac{1}{2h},\quad
\overline{U}_{v_N} ^2 = \frac{1}{4}\,\zeta(v_N)\,
,\quad
\overline{U}_{v_1}\overline{U}_{v_N} =-\overline{U}_{v_N}\overline{U}_{v_1}.
\een
b) The elements
\ben
X:= \Big((-1)^{N-1} \zeta(v_1)\cdots \zeta(v_{N-1})\,
\frac{1}{2h}\Big)^{-1/2} \overline{U}_{v_1} \quad
Y:=\frac{1}{2}\,(-\zeta(v_N))^{-1/2} \overline{U}_{v_N}
\een
generate a subgroup of $G_\sigma$ isomorphic to the quaternion group
$$
Q_8:=\{\pm 1, \pm \mathbf{i}, \pm \mathbf{j},\pm \mathbf{k}\},
$$ 
where the multiplication is given by the standard quaternion
relations $\mathbf{i}^2=\mathbf{j}^2=\mathbf{k}^2=-1$ and $\mathbf{i}\mathbf{j}=\mathbf{k}$.

c) The $\CC^*$-invariant representations of $G_\sigma$ are given by
$\mathbb{H}^n$ ($n\geq 1$), where $\mathbb{H}$ is the quaternion algebra and we fix an identification
$\mathbb{H}\cong \CC^2$, such that the
natural left action of $Q_8$ on $\mathbb{H}$ becomes a 2-dimensional
representation. 
\end{lemma} 
\proof
a) Using that $\sigma^{N-1}(v_1)=-v_1$ we get that 
\ben
(-1)^{N-1} \zeta(v_1)\cdots \zeta(v_{N-1}) U_{v_1}^{-1}U_{-v_1} \in N_\sigma.
\een
On the other hand by definition 
\ben
U_{-v_1}U_{v_1} = \frac{\epsilon(-v_1,v_1)}{B(-v_1,v_1)}\, U_0 =
-\frac{1}{2h} ,
\een
so the relation for $\overline{U}_{v_1}$ follows. The proof of the
second relation is similar, while the last one follows directly from
the definitions.

b) By part a) we have $X^2=-1$, $Y^2=-1$ and $XY=-YX$. We need just to
prove that there are no further relations between $X$ and $Y$. This is
equivalent to proving that elements of the type $c U_{v_1}, cU_{v_N}$,
or $c U_{v_1}U_{v_N}$ do not belong to $N_\sigma$. On the other hand, every element of
$N_\sigma$ has the form $c_\alpha U_{(1-\sigma)\alpha}$ for some
$\alpha\in \ZZ^N$ and $c_\alpha\in \CC^*$. Our claim follows from the
fact that 
\ben
(1-\sigma)^{-1}(v_1)=\frac{1}{2}(v_1+\cdots + v_{N-1}),\quad 
(1-\sigma)^{-1}(v_N)=\frac{1}{2}v_N.
\een

c) Since the action of the Coxeter transformation on $v_1$ and $v_N$
produces two orbits that contain $\pm v_i$ for all $i$, we get that
the group $G_\sigma$ is generated by $\CC^*$, $X$, and $Y$. The
$\CC^*$-invariance implies that every representation is uniquely
determined by its restriction to the quaternion subgroup $Q_8=\langle
X,Y\rangle$. Moreover, the $\CC^*$-invariance implies that $-1\in Q_8$
acts by $-1$ in the representation. The representations of the
quaternion group $Q_8$ are
completely reducible and it is known that the there are 5 irreducible
ones. Their  dimensions are $1,1,1,1$,
and $2$ and only the 2-dimensional one has the property that $-1$ acts by
$-1$.  
\qed

The result of Bakalov and Kac (see \cite{BK}, Proposition 4.2) can be
stated as follows. Suppose that  $M$ is a $\sigma$-twisted
$V_{\ZZ^N}$-module. Note that due to $sigma$-invariance the modes $(H_i)^M_{(-m/h)}$ are not zero only
if $m\equiv m_i (\operatorname{mod} \ h)$. The twisted Borcherd's
identity implies that 
\ben
K\mapsto 1,\quad H_i s^{-n-m_i/h}\mapsto (H_i)^M_{(-n-m_i/h)}\quad
(1\leq i\leq N, n\in \ZZ)
\een 
is a representation of $\widehat{\lieh}_\sigma$. The axioms of a
twisted module imply that 
\ben
Y^M(e^\alpha,\lambda) = U^M_\alpha\, \lambda^{-|\alpha|^2/2} \,
:\exp\Big(\sum_{n\in \frac{1}{h}\ZZ-\{0\}} \alpha^M_{(-n)} \frac{\lambda^n}{n}\Big):
\een
where $U^M_\alpha$ are certain linear operators that commute with the
representation of $\widehat{\lieh}_\sigma$. Note that the
$\sigma$-invariance imply that $U^M_{\sigma\alpha} =\zeta(\alpha)^{-1}
U^M_{\alpha} (-1)^{|\alpha|^2}$, which is the relation imposed on the
elements of the group $G$. Moreover, the assignment
\ben
cU_\alpha\mapsto c U^M_\alpha,\quad c\in \CC^*,\quad \alpha\in \ZZ^N,
\een
defines a representation of $G_\sigma$ on $M$ such that
$U^M_0=\operatorname{Id}_M$. The map just described establishes an
equivalence between the category of $\sigma$-twisted
$V_{\ZZ^N}$-modules and the category of
$(\widehat{\lieh}_\sigma,G_\sigma)$-modules in which $K\in
\widehat{\lieh}_\sigma$ is represented by the
identity operator and the representation of $G_\sigma$ is
$\CC^*$-invariant.

\begin{lemma}\label{le:tw-rep}
The fermionic Fock space $\F$ has a structure of a $\sigma$-twisted
$V_{\ZZ^N}$-module such that 
\begin{align}
\notag
Y^\F(e^{\pm v_i},\lambda) & =  C^{\pm}_i \, \lambda^{-1/2}
\phi_1(\pm\lambda^{1/h} \eta^i),\quad 1\leq i\leq N-1 \\
\notag
Y^\F(e^{\pm v_N},\lambda) & = C^{\pm}_N\, \lambda^{-1/2}
\phi_2(\pm \lambda^{1/2}),
\end{align}
where $C^\pm_i$ ($1\leq i\leq N$) are some constants satisfying 
\ben
C^+_i\, C^-_i= -1/h\quad  (1\leq i\leq N-1),\quad C_N^+\, C_N^-=-1/2.
\een
\end{lemma}
\proof
Recalling the commutation relations for $J^a_m$ we get that the map
\ben
H_i s^{-m/h}\mapsto \frac{1}{\sqrt{2}} J^1_{-m}\quad (1\leq i\leq N-1)\quad
H_N s^{-m(N-1)/h}\mapsto \sqrt{\frac{N-1}{2}} J^2_{-m}
\een
defines a representation of the twisted Heisenberg algebra
$\widehat{\lieh}_\sigma$ in which $K\mapsto 1$. Since the operators $Q_1$ and $Q_2$ satisfy the relations
$Q_1^2=Q_2^2=1/2$ and $Q_1Q_2=-Q_2Q_1$, we can choose constants
$C_1^+$ and $C_N^+$ such that the map 
\begin{align}
\notag
X& \mapsto \Big((-1)^{N-1}\zeta(v_1)\cdots \zeta(v_{N-1})\,\frac{1}{2h}\Big)^{-1/2}
C_1^+ Q_1,\\
\notag
Y & \mapsto \frac{1}{2} \Big( -\zeta(v_N)\Big)^{-1/2} C_N^+ Q_2
\end{align}
defines a representation of the quaternion group $Q_8$ and hence a
$\CC^*$-invariant representation of $G_\sigma$ on the
fermionic Fock space $\mathcal{F}$.  Therefore, we can equip
$\mathcal{F}$ with the structure of a $\sigma$-twisted
$V_{\ZZ^N}$-module, such that $Y^{\mathcal{F}}(e^{\pm v_i},\lambda)$ has the
form stated in the lemma. 

Let us porve the relations between the constants. Note that $U^\F_{v_1}= C_1^+ Q_1$ and
$U^\F_{-v_1}=C_1^-Q_1$, so 
\ben
-\frac{1}{2h}= U^\F_{v_1} U^\F_{-v_1} = C_1^+ C_1^-\, \frac{1}{2}. 
\een
The proof of the remaining relations is similar.
\qed

\subsection{Proof of Theorem \ref{t1}}
The Casimir $\Omega_{\rm KW}$ of the Kac--Wakimoto hierarchy is given by the residue of
the 1-form $Y^\F(\omega_{\rm KW},\lambda)\lambda d\lambda$. 
Recall Lemma \ref{le:casimir} and the product formula
\eqref{ope-voa}. After a direct computation using Lemma
\ref{le:tw-rep} we get that 
\ben
Y^\F(\omega_{\rm BKP},\lambda) = \sum_{m\in \ZZ}
\widetilde{\Omega}_m \lambda^{-m-1},
\een
where 
\ben
\widetilde{\Omega}_m =-\sum_{a=1}^2 \sum_{k\in \ZZ} (-1)^k
\phi_a(k)\otimes \phi_a(-k-m h_a).
\een
Moreover another straightforward computation using the identity
\ben
[Y^\F(a,\lambda_1), Y^\F(b,\lambda_2) ]= \sum_{n=0}^\infty
\frac{1}{n!}
Y^\F(a_{(n)} b,\lambda_2)  \partial_{\lambda_2}^n \delta(\lambda_1,\lambda_2),
\een
where $a,b\in V_{\ZZ^N}$ are $\sigma$-invariant and
$\delta(\lambda_1,\lambda_2):=\sum_{n\in \ZZ} \lambda_1^n 
\lambda_2^{-n-1}$ is the formal delta function,  yields the following
commutation relations
\ben
[\widetilde{\Omega}_m , \widetilde{\Omega}_n]=2N m \delta_{m,-n}.
\een
Let us define normal ordering 
\ben
:\widetilde{\Omega}_m \widetilde{\Omega}_n: =
\begin{cases}
\widetilde{\Omega}_m \widetilde{\Omega}_n & \mbox{ if $n\geq 0$}, \\
\widetilde{\Omega}_n \widetilde{\Omega}_m & \mbox{ if $n<0$}.
\end{cases}
\een
Then we have
\ben
Y^\F(\omega_{\rm BKP},\lambda_1) 
Y^\F(\omega_{\rm BKP},\lambda) = :  
Y^\F(\omega_{\rm BKP},\lambda_1) 
Y^\F(\omega_{\rm BKP},\lambda): + \frac{2N}{(\lambda_1-\lambda_2)^2}.
\een
Recalling the product formula (with $k=2$) we get 
\ben
-\frac{1}{2}Y^\F((\omega_{\rm BKP})_{(-1)} \omega_{\rm
  BKP},\lambda)= 
-\frac{1}{2} :  
Y^\F(\omega_{\rm BKP},\lambda) 
Y^\F(\omega_{\rm BKP},\lambda):\ .
\een
We get the following formula for the Casimir
\ben
\Omega_{\rm KW} = -\frac{1}{2} \, \widetilde{\Omega}_0^2
-\sum_{m=1}^\infty \widetilde{\Omega}_{-m}\widetilde{\Omega}_m.
\een
Finally, it remains only to recall that by definition 
\ben
\widetilde{\Omega}_m = -(Q_1\otimes Q_1) \, \Omega_m
\een
and that $(Q_1\otimes Q_1)^2=\frac{1}{4}$. \qed

\section{Virasoro symmetries}

The goal in this section is to prove Corollary \ref{c1} and Theorem
\ref{t2}. We will need to work out the commutation relations between
the fermions $\phi_a(k)$ and the Virasoro
operators $D_k$ (see \eqref{vir-D_k}). This can be done directly (see
for example \cite{KvL}). For the sake of completeness we would like
to use the twisted module structure on $\F$ and derive all commutation
relations from the Borcherd's identities.

\subsection{Commutation relations}

Recall the $\sigma$-twisted $V_{\ZZ^N}$-module structure on $\F$ from
the previous section (see Lemma \ref{le:tw-rep}). By definition we
have
\ben
Y^\F(H_N,\lambda)=\sqrt{\frac{N-1}{2}} \sum_{m\in \ZZ_{\rm odd}} J^2_{-m} \lambda^{m/2-1}
\een
and 
\ben
Y^\F(H_i,\lambda) = \frac{1}{\sqrt{2}} \sum_{m} J^{1}_{-m}
\lambda^{m/h-1},\quad 1\leq i\leq N-1,
\een 
where the sum is over all $m\in \ZZ$ such that $m\equiv
m_i(\operatorname{mod}\ h)$. Therefore, the twisted fields
representing the standard basis $v_i$ are given by 
\ben
Y^\F(v_i,\lambda) = \frac{1}{h\lambda} \sum_{m\in \ZZ_{\rm odd}}
J^1_{-m} (\lambda^{1/h}\eta^i)^m,\quad 1\leq i\leq N-1,
\een
and 
\ben
Y^\F(v_N,\lambda) = \frac{1}{2\lambda} \sum_{m\in \ZZ_{\rm odd}}
J^2_{-m} (\lambda^{1/2})^m.
\een
\begin{lemma}
The conformal vector \eqref{conf-vector} is represented by the twisted field
\ben
Y^\F(\nu,\lambda) =\sum_{k\in \ZZ} D_k\lambda^{-k-2},
\een
where $D_k$ are the operators \eqref{vir-D_k}.
\end{lemma}
\proof
This is a straightforward computation using the product formula
\eqref{ope-voa}. Let us point out the main steps leaving some of the details
to the reader. We have 
\ben
&&
Y^\F(v_i,\lambda_1) Y^\F(v_i,\lambda_2) =\\
&&
 :Y^\F(v_i,\lambda_1)Y(v_i,\lambda_2):
+ 
\frac{\lambda_1^{1/h-1}\lambda_2^{1/h-1}  }{(\lambda_1^{1/h}-\lambda_2^{1/h})^2h^2} +
\frac{\lambda_1^{1/h-1}\lambda_2^{1/h-1} }{(\lambda_1^{1/h}+\lambda_2^{1/h})^2h^2},
\een 
where $1\leq i\leq N-1$ and the normal ordering means that the currents $J^1_m$ with $m>0$
should be applied first. The formula for $i=N$ is the same except that
we have to replace everywhere $h$ with $2$.

Note that we have the following Taylor's series expansions at
$\lambda_1=\lambda_2$:
\ben
(\lambda_1-\lambda_2)^2\, 
\frac{ \lambda_1^{1/h-1}\lambda_2^{1/h-1}
}{(\lambda_1^{1/h}-\lambda_2^{1/h})^2h^2} =
1+ \frac{h^2-1}{12 h^2\lambda_2^2}\, (\lambda_1-\lambda_2)^2 + \cdots 
\een
and 
\ben
(\lambda_1-\lambda_2)^2\, 
\frac{ \lambda_1^{1/h-1}\lambda_2^{1/h-1}
}{(\lambda_1^{1/h}+\lambda_2^{1/h})^2h^2} =
\frac{1}{4h^2\lambda_2^2}\, (\lambda_1-\lambda_2)^2 + \cdots\ .
\een
Recalling the product formula \eqref{ope-voa} we get 
\ben
Y^\F(v_{i(-1)}v_i,\lambda) = :Y^\F(v_i,\lambda) Y^\F(v_i,\lambda): \
+\  
\frac{h^2+2}{12h^2} \, \lambda^{-2}
\een
for all $1\leq i\leq N-1$. For $i=N$ the formula remains the same
except that we have to replace $h$ with $2$. 
Recalling the formulas for $Y^\F(v_i,\lambda)$ we get
\ben
\frac{1}{2} \sum_{i=1}^{N-1} :Y^\F(v_i,\lambda) Y^\F(v_i,\lambda): =
\sum_{k\in \ZZ} \sum_{m\in \ZZ_{\rm odd}} 
\frac{1}{4h} :J^1_m J^1_{-m-kh}: \lambda^{-k-2},
\een
and
\ben
\frac{1}{2} 
:Y^\F(v_N,\lambda) Y^\F(v_N,\lambda): =
\sum_{k\in \ZZ} \sum_{m\in \ZZ_{\rm odd}} 
\frac{1}{8} :J^2_m J^2_{-m-2k}: \lambda^{-k-2}.
\een
Finally it remains only to sum up the extra contributions to the
coefficient in front of $\lambda^{-2}$: 
\ben
\frac{h^2+2}{24 h^2}\, (N-1) + \frac{2^2+2}{24\cdot 2^2} = \frac{(h+1)N}{24h}.\qed
\een
The twisted Borcherd's identity \eqref{tw_borcherd-id} with $n=0$
imply the following formula
\ben
[Y^\F(\nu,\lambda_1),Y^\F(a,\lambda_2)] = \sum_{m=0}^\infty
\frac{1}{m!} Y^\F(\nu_{(m)}a,\lambda_2)\, 
\partial_{\lambda_2}^m\delta(\lambda_1,\lambda_2),
\een
where $a\in V_{\ZZ^N}$ and we used that $\sigma(\nu)=\nu.$ 
\begin{lemma}\label{le:vir-fermion_comm}
The following formulas hold:
\ben
[D_k,\phi_a(l)]= \Big(\frac{l}{h_a} -\frac{k}{2}\Big) \phi_a(l-kh_a)
\een
and
\ben
[D_k,\phi_a(z_a)]= 
\Big(\frac{1}{h_a} z_a^{1+kh_a}\partial_{z_a}
+\frac{k}{2}\, z_a^{kh_a} \Big) \phi_a(z_a),\quad a=1,2.
\een
\end{lemma}
\proof
The first formula is equivalent to the second one as one can see
immediately by comparing the coefficients in front of $z_a^l$. Let us
prove the second formula. We will consider only the case $a=1$. The
argument for $a=2$ is similar. 

Since 
\ben
\nu_{(0)}e^{v_1} = v_{1(-1)} e^{v_1},\quad \nu_{(1)} e^{v_1} =
\frac{1}{2} e^{v_1},\quad \nu_{(m)} e^{v_1}=0 \ (m>1),
\een
we get 
\ben
[Y^\F(\nu,\lambda_1),Y^\F(e^{v_1},\lambda_2)] = 
Y^\F(v_{1(-1)}e^{v_1},\lambda_2)\delta(\lambda_1,\lambda_2)
+\frac{1}{2}
Y^\F(e^{v_1},\lambda_2) \partial_{\lambda_2}\delta(\lambda_1,\lambda_2). 
\een
After a straightforward computation we get 
\ben
&&
Y^\F(v_1,\lambda_1)Y^\F(e^{v_1},\lambda_2) = \\
&&
:Y^\F(v_1,\lambda_1)Y^\F(e^{v_1},\lambda_2) :
+\Big(\frac{\lambda_1^{1/h-1}}{(\lambda_1^{1/h}-\lambda_2^{1/h})h} -
\frac{\lambda_1^{1/h-1}}{(\lambda_1^{1/h}+\lambda_2^{1/h})h} 
\Big)\, Y^\F(e^{v_1},\lambda_2) .
\een
Using the Taylor's series expansion 
\ben
(\lambda_1-\lambda_2)
\Big(
\frac{\lambda_1^{1/h-1}}{(\lambda_1^{1/h}-\lambda_2^{1/h})h} -
\frac{\lambda_1^{1/h-1}}{(\lambda_1^{1/h}+\lambda_2^{1/h})h} \Big)=
1 -\frac{1}{2\lambda_2} \, (\lambda_1-\lambda_2) +\cdots 
\een
and the product formula \eqref{ope-voa} we get 
\ben
Y^\F(v_{1(-1)}e^{v_1},\lambda)= \partial_\lambda \, Y^\F(e^{v_1},\lambda),
\een
where we used that 
\ben
:Y^\F(v_1,\lambda)Y^\F(e^{v_1},\lambda) : =
\Big( \partial_\lambda
+\frac{1}{2\lambda}\Big)\,
Y^\F(e^{v_1},\lambda).
\een
The formula for the commutator that we want to compute takes the form
\ben
[Y^\F(\nu,\lambda_1), Y^\F(e^{v_1},\lambda_2)] = 
\Big( \delta(\lambda_1,\lambda_2)
\partial_{\lambda_2}+
\frac{1}{2}\partial_{\lambda_2}\delta(\lambda_1,\lambda_2) 
\Big)\, Y^\F(e^{v_1},\lambda_2). 
\een
Comparing the coefficients in front of $\lambda_1^{-k-2}$ we get 
\ben
[D_k, Y^\F(e^{v_1},\lambda_2)]=
\Big( \lambda_2^{k+1}\partial_{\lambda_2}+
\frac{(k+1)}{2}\, \lambda_2^k\Big)\, Y^\F(e^{v_1},\lambda_2).
\een
It remains only to recall that $Y^\F(e^{v_1},\lambda_2) = C_1^+
\lambda_2^{-1/2} \phi_1(\lambda_2^{1/h} \eta)$ and to substitute
$z_1=\lambda_2^{1/h}\eta$.
\qed 

\subsection{Proof of Corollary \ref{c1}}

Let us prove part a) of Corollary \ref{c1}.
We have to prove that if $\Omega_{\rm KW}(\tau\otimes \tau)=0$ then
$\Omega_m(\tau\otimes\tau)=0$ for all $m\geq 0$. Following ten
Kroode--van de Leur (see \cite{KvL}) we equip the Fock space $\F$
with a positive definite Hermitian form $H$, such that 
\begin{enumerate}
\item
$H(|0\rangle,|0\rangle)=1$
\item
$H(\phi_a(k) v_1,v_2)=(-1)^k H(v_1,\phi_a(-k)v_2)$ for all $k\in \ZZ$,
$a=1,2$, and $v_1,v_2\in \F$. 
\end{enumerate}
In particual $\F^{\otimes 2}$ also has an induced positive definite
Hermitian form. Observe that if $\Omega_{\rm KW}(v)=0$ for some $v\in
\F^{\otimes 2}$, 
then since the Hermitian adjoint of the operator $\Omega_{-m}$ is $\Omega_m$ we
get
\ben
-4H(\Omega_{\rm KW}(v), v) =
\frac{1}{2}H(\Omega_{0}(v),
\Omega_{0}(v))+
\sum_{m=1}^\infty  
H(\Omega_{m}(v),
\Omega_{m}(v)).
\een
The positive definitness of $H$ implies that
$\Omega_m(v)=0$ for all $m\geq 0$. 

We can not apply the above argument directly, because $\tau$ is a
formal power series. It belongs to the completion
$\widehat{\F_0}=\CC[\![\mathbf{t}^1,\mathbf{t}^2]\!]$ and the
Hermitian form does not extend to $\widehat{\F}_0$. On the other hand
the grading operator
\ben
D_0-\frac{(h+1)N}{24h} = \sum_{m\in \ZZ_{\rm odd}^+}
\frac{m}{h} t^1_m \partial_{t^1_m} + 
\frac{m}{2} t^2_m \partial_{t^2_m} 
\een
commutes with $\Omega_{\rm KW}$ (see Lemma
\ref{le:vir-fermion_comm}) and $\widehat{\F}_0$ decomposes as an
infinite product of finite dimensional eigensubspaces. Therefore if we decompose
$\tau\otimes \tau=\sum_{n\geq 0} v_n$ where each $v_n\in \F_0$ is
homogeneous of degree $n$, then $\Omega_{\rm KW}(v_n)=0$ for all
$n$. Recalling the argument from above $\Omega_m(v_n)=0$ for all
$m\geq 0$, i.e., $\Omega_m(\tau\otimes \tau)=0.$   

Let us prove part b). Suppose that $\tau$ is a tau-function of the
Kac--Wakimoto hierarchy. By part a) $\Omega_m(\tau\otimes \tau)=0$ for
all $m\geq 0$. In particular $\tau$ is a tau-function of the
2-component BKP. Let us define the wave function 
\ben
\Psi(\mathbf{t},z)=\Psi^{(1)}(\mathbf{t},z_1)e_1+\mathbf{i} \Psi^{(2)}(\mathbf{t},z_2)e_2,
\een
where
\ben
\Psi^{(a)}(\mathbf{t},z_a)=\Gamma(\mathbf{t}^a,z_a)\tau/\tau,\quad a=1,2.
\een
The subspace $U\in \operatorname{Gr}_2^{I,(0)}$ is spanned by the
coefficients of the Taylor's series expansion of $\Psi(\mathbf{t},z)$
at $\mathbf{t}=0$. Therefore we need to prove that 
$(z_1^h,z_2^2) \Psi(\mathbf{t},z)\in U$. 

By definition $\Omega_m$ is the following bi-linear operator acting on
$\CC[\![\mathbf{t}^1,\mathbf{t}^2]\!] ^{\otimes 2}$
\ben
\operatorname{Res}_{z_1=0} \frac{dz_1}{z_1} z_1^{mh}\Big(
  \Gamma(\mathbf{t}^1,z_1) \otimes \Gamma(\mathbf{t}^1,-z_1) \Big)
  -
\operatorname{Res}_{z_2=0} \frac{dz_2}{z_2} z_2^{2m}\Big(
  \Gamma(\mathbf{t}^2,z_2) \otimes \Gamma(\mathbf{t}^2,-z_2) \Big).
\een
The equation $\Omega_m(\tau\otimes \tau)=0$  is equivalent to 
\ben
((z_1^{mh},z_2^{2m})\Psi(\mathbf{t}',z), \Psi(\mathbf{t}'',z) =0.
\een
Therefore $(z_1^{mh},z_2^{2m})\Psi(\mathbf{t}',z)$ is orthogonal to
every vector in $U$. Since $U$ is maximally isotropic, we must have
$(z_1^{mh},z_2^{2m})\Psi(\mathbf{t}',z)\in U$. Note that the argument
is invertible, i.e., the condition that $(z_1^h,z_2^2)U\subset U$
implies that the tau-function  $\tau$ satisfies all bi-linear
equations $\Omega_m(\tau\otimes \tau)=0$, so by Theorem \ref{t1}
$\tau$ is a tau-function of the Kac--Wakimoto hierarchy.

\subsection{Virasoro constraints}
Suppose now that $\tau$ is a tau-function of the Kac--Wakimoto
hierarchy satisfying the Virasoro constraints $L_k\tau=0$ for all
$k\geq -1$, where $L_k$ are the operators \eqref{vir-L_k}.
Recall the differential operators $\ell_k(z)$ (see formula \eqref{vir-ell}).
\begin{lemma}\label{vir-Gr_constr}
Let $U\in \operatorname{Gr}_2^{I,(0)}$ be the point corresponding to
a tau-function $\tau$ of the Kac--Wakimoto hierarchy. Then the condition $L_k(\tau)=0$ implies
$\ell_k(z) U\subset U$.
\end{lemma}
\proof
Note that 
\ben
[J^1_{1+(1+k)h},\Gamma(\mathbf{t}^a,z_a)]=2\delta_{1,a} z_1^{1+(1+k)h}\Gamma(\mathbf{t}^a,z_a).
\een
Recalling the commutation relations in Lemma \ref{le:vir-fermion_comm}
and the boson-fermion isomorphism $\phi_a(z_a)=Q_a
\Gamma(\mathbf{t}^a,z_a)$ we get
\ben
[L_k, \Gamma(\mathbf{t}^a,z_a)]=\ell_k^{(a)}(z_a) \, \Gamma(\mathbf{t}^a,z_a).
\een
If $\Psi(\mathbf{t},z)$ is the wave function then 
\ben
\ell_k(z) \Psi(\mathbf{t},z) = 
\tau^{-1}(\mathbf{t})\, \ell_k(z) \, (\Gamma(\mathbf{t}^1,z_1) e_1
+\mathbf{i}  \Gamma(\mathbf{t}^2,z_2) e_2)\, \tau(\mathbf{t}).
\een
Since 
\ben
&&
\ell_k(z) \, (\Gamma(\mathbf{t}^1,z_1) e_1
+\mathbf{i}  \Gamma(\mathbf{t}^2,z_2) e_2) = \\
&&
L_k\circ (\Gamma(\mathbf{t}^1,z_1) e_1
+\mathbf{i}  \Gamma(\mathbf{t}^2,z_2) e_2)-
(\Gamma(\mathbf{t}^1,z_1) e_1
+\mathbf{i}  \Gamma(\mathbf{t}^2,z_2) e_2)\circ L_k
\een
and $L_k$ anihilates $\tau$ we get
\beq\label{ell_k-constr}
\ell_k(z) \Psi(\mathbf{t},z) = (\tau^{-1} \circ L_k \circ \tau) \Psi(\mathbf{t},z).
\eeq
The RHS of the above formula belongs to $U$, because
$\Psi(\mathbf{t},z)$ does and $\tau^{-1} \circ L_k \circ \tau$ is a
differential operator in the dynamical variables $\mathbf{t}$. 
\qed
\begin{lemma}\label{le:string-p}
If $p\geq 0$ is an integer, then the following formulas hold:
\ben
(\mathbf{i}\ell^{(1)}_{-1}(z_1) )^p \Psi^{(1)}(x,z_1) = z_1^p
\Psi^{(1)}(x,z_1)+ \Big(\frac{\mathbf{i}p}{h}\, x_1 z_1^{p-h} +
O(z_1^{p-h-1})\Big) e^{x_1 z_1}
\een
and 
\ben
(\mathbf{i}\ell^{(2)}_{-1}(z_2) )^p \Psi^{(2)}(x,z_2) = O(z_2^{-p})
e^{x_2 z_2}.
\een
\end{lemma}
\proof
The second formula is obvious. Let us prove the first one. We argue by
induction on $p$. For $p=0$ the identity is obvious. 
Suppose the formula is true for $p$. We will prove it for $p+1$. 
To begin with note that
\ben
\ell_{-1}^{(1)}(z_1)= 
z_1 -\frac{\mathbf{i}}{2} z_1^{-h}+\frac{\mathbf{i}}{h} z_1^{-h} (z_1\partial_{z_1}).
\een
We have
\ben
&&
\Big( z_1 -\frac{\mathbf{i}}{2} z_1^{-h}+\frac{\mathbf{i}}{h} z_1^{-h}
(z_1\partial_{z_1})\Big) z_1^p \Psi^{(1)}(x,z_1) =\\
&&
z_1^p\Big( z_1 \, \Psi^{(1)}(x,z_1) +
\Big(\frac{\mathbf{i}}{h} x_1 z_1^{1-h}+O(z_1^{-h})\Big) e^{x_1
  z_1}\Big) + O(z_1^{p-h}) e^{x_1z_1}=\\
&&
z_1^{p+1} \Psi^{(1)}(x,z_1) +\Big(\frac{\mathbf{i}}{h} x_1 z_1^{p+1-h}
+ O(z_1^{p-h})\Big) e^{x_1z_1},
\een
where we used that $\Psi^{(1)}(x,z_1)=(1 +O(z_1^{-1})) e^{x_1 z_1}$.
We also have
\ben
&&
\Big( z_1 -\frac{\mathbf{i}}{2} z_1^{-h}+\frac{\mathbf{i}}{h} z_1^{-h}
(z_1\partial_{z_1})\Big) \Big(
\frac{\mathbf{i}p}{h} x_1 z_1^{p-h}\, e^{x_1z_1}\Big) = \\
&&
\Big(
\frac{\mathbf{i}p}{h} x_1 z_1^{p+1-h} + O(z_1^{p-2h+1})\Big)\, e^{x_1z_1}.
\een
Note that $p-2h+1\leq p-h$, so combining the above two formulas
completes the inductive step.
\qed

Let us prove the first equation in Theorem \ref{t2}. We will need the
following simple lemma.
\begin{lemma}\label{le:int}
Suppose that $U\in \operatorname{Gr}_2^{I,(0)}$ and that the series 
\ben
w(x,z)=(u_1(x,z_1) e^{x_1 z_1},u_2(x,z_2) e^{x_2 z_2})
=\sum_{m,n=0}^\infty w_{m,n}(z) x_1^m x_2^n
\een
belongs to the intersection
\ben
U\cap (e^{x_1z_1},e^{x_2z_2}) \cdot V_0,
\een
where the multiplication $\cdot$ is the componentwise multiplication. 
Then $w=0$.
\end{lemma}
\proof
Using Taylor's formula we get
\ben
w_{m,n}=\left.\frac{1}{m!n!} ((\partial_1+z_1)^m\partial_2^n
u_1, \partial_1^m (\partial_2+z_2)^n u_2)\right|_{x_1=x_2=0}.
\een 
Using the binomial formula we get
\ben
w_{m,n}=
\sum_{i=0}^m u_{1,i,n}(z_1) \frac{z_1^{m-i}}{(m-i)!} e_1 + 
\sum_{j=0}^n u_{2,m,j}(z_2) \frac{z_2^{n-j}}{(n-j)!} e_2,
\een
where $u_{a,m,n}(z_a)$ is the coefficient in front of $x_1^m x_2^n$ of
$u_a(x,z_a)$.
We argue by induction on $m+n$ that
$u_{a,m,n}=0$ for all $a=1,2$ and $m,n\geq 0.$ If $m=n=0$ the definition of the
Grassmanian implies that 
$
(u_{1,0,0},u_{2,0,0})=w_{0,0}\in U\cap V_0=\{0\}.
$
Suppose $u_{a,m,n}=0$ for $m+n<k$. If $m+n=k$ then  $u_{1,i,n}=0$ and
$u_{2,m,j}=0$ for $i<m$ and $j<n$, respectively. Hence 
\ben
w_{m,n}=(u_{1,m,n},u_{2,m,n}) \in U\cap V_0=\{0\},
\een
where we used that $(u_1(x,z_1),u_2(x,z_2))\in V_0$.
\qed

\begin{lemma}\label{tau-initial}
If $\tau(\mathbf{t})$ satisfies the Virasoro and the dilaton constraints
then the restriction of $\tau$ to $t^a_m=0$ for all $a=1,2$ and $m>1$
is given by $e^{c-\mathbf{i} x_1 x_2^2/8}$, where $c$ is some constant
independent of $x_1$ and $x_2$.
\end{lemma}
\proof
The restriction of the string equation to $t^a_m=0$ for $a=1,2$, $m>1$
yields 
\ben
-\mathbf{i}\partial_1\tau(x_1,x_2) + \frac{1}{8} x_2^2 \tau(x_1,x_2)=0.
\een
Therefore $\tau(x) = e^{f(x_2)-\mathbf{i} x_1x_2^2/8 }$ for some
function $f(x_2)\in \CC[\![x_2]\!]$. We have to check that $f$ is a constant. 
Let us write $\tau(\mathbf{t})=\mathcal{D}(1,\mathbf{t})$ where
$\mathcal{D}(\hbar,\mathbf{t})$ is a solution of the dilaton equation
\eqref{dil-eq}. The differential operator defining the dilaton
equation commutes with 
$$
L_0= -\mathbf{i}\partial_{t^1_{1+h}} + \sum_{m\in \ZZ^+_{\rm odd}}
(\frac{m}{h} t^1_m \partial_{t^1_m} +
\frac{m}{2}t^2_m\partial_{t^2_m}) + \frac{N(h+1)}{24h}. 
$$
Therefore
$L_0\mathcal{D}=0$. Subtracting the dilaton equation from
$\frac{h}{h+1} L_0\mathcal{D}=0$ we get 
\ben
\Big(
-2\hbar\partial_\hbar +
\frac{1}{h+1}\, \sum_{m\in \ZZ^+_{\rm odd}} 
((m-1-h)t^1_m\partial_{t^1_m} + ((N-1) m-1-h) t^2_m \partial_{t^2_m}
\Big) \mathcal{D}(\hbar,\mathbf{t})=0. 
\een
This differential equation imposes a non-trivial constraint on the
coefficients of $\mathcal{D} $. Namely, let us write 
\ben
\mathcal{D}(\hbar,\mathbf{t}) = \exp\Big(
\sum_{g,\mu,\nu} C^{(g)}(\mu,\nu) 
t^1_{m_1}\cdots t^1_{m_r} 
t^2_{n_1}\cdots t^2_{n_s} \, \hbar^{g-1}
\Big)
\een
where the sum is over all $g\geq 0$, $\mu=(m_1,m_2,\dots,m_r)$, and
$\nu=(n_1,n_2,\dots,n_s)$.  Then if $C^{(g)}(\mu,\nu)\neq 0$ we must
have
\ben
\frac{1}{h+1} \Big( \sum_{i=1}^r m_i +\sum_{j=1}^s (N-1) n_j\Big)
=2g-2 + r + s.
\een
In order to prove that $f(x_2)$ is a constant independent of $x_2$ we
just need to check that we can 
not have monomials with $r=0$ and $n_j=1$ for $1\leq j\leq s$. Suppose
that such a monomial contributes to the tau-function. Then 
\ben
\frac{s(N-1)}{2N-1} = 2g-2 +s.
\een
Since $N-1$ and $2N-1$ are relatively prime, we get $s=\ell (2N-1)$
for some integer $\ell\geq 1$. The above equation takes the form 
$\ell (N-1)=2g-2+\ell(2N-1)$, i.e., $2g-2+\ell N =0$. This equation
does not have solutions for $g\geq 0$ and $N\geq 4$. Hence $\log
\tau(x)$ does not have non-trivial monomials in $x_2$.
\qed 

\smallskip

For future references let us state an important corollary of the proof
of Lemma \ref{tau-initial}. We got that if some monomial 
$
t^1_{m_1}\cdots t^1_{m_r} 
t^2_{n_1}\cdots t^2_{n_s}
$
contributes to the tau-function, then the number $\sum_{i=1}^r m_i
+\sum_{j=1}^s (N-1) n_j$ is an integer divisible by $h+1$. In
particular, we get that the total descendant potential is invariant
under the rescaling 
\beq\label{scaling-(h+1)}
t^1_m \mapsto \omega^m\, t^1_m,\quad t^2_m\mapsto \omega^{m(N-1)}
t^2_m \ ,
\eeq
where $\omega=e^{2\pi\mathbf{i}/(h+1)}$.

\subsection{Proof of Theorem \ref{t2}} 
The third equation is a direct consequence of Lemma
\ref{tau-initial}. Indeed, by property (W2) of
wave functions we know that there exists $q(x)\in \CC[\![x_1,x_2]\!]$
such that $(\partial_1\partial_2+q(x)) \Psi(x,z)=0,$ where recall that
$x_1:=t^1_1$, $x_2:=t^2_1$, and $\partial_a=\partial/\partial
x_a.$. By definition   
\ben
\Psi^{(1)}(x,z_1) = 
\Big(1-2\frac{\partial_1\tau}{\tau} z_1^{-1} +\cdots \Big)
e^{x_1 z_1}.
\een
Comparing the coeffiecients in front of $z_1^0$ in 
$(\partial_1\partial_2+q(x)) \Psi^{(1)}(x,z_1)=0$ we get 
$$
q(x)=2\partial_1\partial_2\log \tau = -\frac{\mathbf{i}}{2}\, x_2.
$$
Let us prove the second equation. Note that the operator $L_{-1}$ is given explicitly by 
\begin{align}
\notag
L_{-1} = & -\mathbf{i}\partial_{t_1^1} + \sum_{a=1,2}\sum_{2i+2j=h_a+2}\frac{(2i-1)(2j-1)}{4h_a}\,
t^a_{2i-1}t^a_{2j-1} +\\
\notag
&
\sum_{a=1,2}\frac{1}{h_a}\sum_{m\in \ZZ_{\rm odd}^+}(m+h_a)t^a_{m+h_a} \partial_{t^a_m}.
\end{align}
By using Lemma \ref{tau-initial} and restricting the equation \eqref{ell_k-constr} with $k=-1$ to $t^a_m=0$
for $a=1,2$ and $m>1$, we get the second equation in
Theorem \ref{t2}:
\beq\label{eq-2}
\ell_{-1}(z)\Psi(x,z) = -\mathbf{i} \partial_1 \Psi(x,z).
\eeq
It remains to prove the first equation. Recalling Lemma
\ref{le:string-p} with $p=h$ we get
\ben
(\partial_1^h +\partial_2^2 -\mathbf{i} x_1 -z_1^h)\Psi^{(1)}(x,z_1) =
O(z_1^{-1}) e^{x_1 z_1}.
\een 
By definition 
\ben
\Psi^{(2)}(x,z_2) = 
\Big(1-2\frac{\partial_2\tau}{\tau} z_2^{-1} +\cdots \Big)
e^{x_2 z_2} = \Big(1+\frac{\mathbf{i}}{2} x_1 x_2 z_2^{-1} +\cdots \Big)
e^{x_2 z_2} ,
\een
where in the second equality we used Lemma \ref{tau-initial}. It is
straightforward to check that 
\ben
(\partial_1^h +\partial_2^2 -\mathbf{i} x_1 -z_2^2)\Psi^{(2)}(x,z_2) =
O(z_2^{-1}) e^{x_2 z_2}.
\een
Hence
\ben
(\partial_1^h +\partial_2^2 -\mathbf{i} x_1 -(z_1^h,z_2^2))\Psi(x,z)
\in U\cap (e^{x_1 z_1},e^{x_2 z_2}) V_0.
\een
It remains only to recall Lemma \ref{le:int}.
\qed

\subsection{The wave function}

Our goal is to prove Corollary \ref{c2}. 
We are going to prove that the system of differntial equations in
Theorem \ref{t2} uniquely determines the wave function $\Psi(x,z)$. 

Using the method of the characteristics it is straightforward to check
that the components of a wave function solving the linear PDEs
\ben
\partial_1 \Psi^{(a)}(x,z_a) = \mathbf{i} \ell_{-1}^{(a)} (z_a)
\Psi(x,z_a),\quad a=1,2,
\een
must have the following form
\ben
\Psi^{(1)}(x,z_1) = e^{
\frac{
\mathbf{i}h
}{
h+1} 
(z_1^{h+1} -(\mathbf{i}  x_1+z_1^h)^{\frac{h+1}{h}} )
}
\frac{z_1^{h/2}}{(\mathbf{i} x_1+z_1^h)^{1/2}}
\Big( 
1 + 
\sum_{k=1}^\infty 
\frac{\psi_k^{(1)}(x_2)}{(\mathbf{i} x_1 +  z_1^h)^{k/h}}
\Big) 
\een
and 
\ben
\Psi^{(2)}(x,z_2) = e^{x_2 (\mathbf{i}  x_1+z_2^2)^{1/2}} 
\frac{z_2}{(\mathbf{i} x_1+z_2^2)^{1/2}}
\Big( 
1 + 
\sum_{k=1}^\infty 
\frac{\psi_k^{(2)}(x_2)}{(\mathbf{i} x_1 +  z_2^2)^{k/2}}
\Big) ,
\een
where $\psi_k^{(a)}(x_2)\in \CC[\![x_2]\!]$ are to be determined from
the remaining two equations.

Let us prove that the 3rd equation uniquely determines $\Psi(x,z)$
from its restriction to $x_2=0$. Conjugating the differential operator
$\partial_1\partial_2-\frac{\mathbf{i}}{2} x_2$ with the exponential
factor of $\Psi^{(1)}$, we get that the series 
\ben
\sum_{k=0}^\infty \psi^{(1)}_k(x_2) (\mathbf{i} x_1 +
z_1^h)^{-k/h-1/2},\quad \psi^{(1)}_0(x_2):=1,
\een
is anihilated by the differential operator
\ben
(\partial_1 + (\mathbf{i} x_1 + z_1^h)^{1/h}) \partial_2 - \frac{\mathbf{i}}{2} x_2.
\een
Comparing the coefficients in front of 
$(\mathbf{i} x_1 + z_1^h)^{-k/h-1/2}$, we get the following recursion
relation
\beq\label{rec-reln-1}
\mathbf{i}\Big(-\frac{k}{h}+\frac{1}{2}\Big) \partial_2\psi^{(1)}_{k-h}(x_2)
+\partial_2\psi^{(1)}_{k+1}(x_2) -\frac{\mathbf{i}}{2} x_2\,
\psi^{(1)}_{k}(x_2) =0.
\eeq
Assuming that we have determined $\psi^{(1)}_{k}(0)$ for all $k$, then
the relation \eqref{rec-reln-1} allows us to reconstruct recursively
$\psi^{(1)}_{k}(x_2)$ for all $k$. 
Similarly the series 
\ben
\sum_{k=0}^\infty \psi^{(2)}_k(x_2) (\mathbf{i} x_1 +
z_2^2)^{-k/2-1/2},\quad \psi^{(2)}_0(x_2):=1,
\een
is anihilated by the differential operator
\ben
&&
\Big(\partial_1+\frac{\mathbf{i}}{2} x_2 (\mathbf{i} x_1 +
z_2^2)^{-1/2}\Big)
\Big(\partial_2 + (\mathbf{i} x_1 + z_2^2)^{1/2} \Big)
-\frac{\mathbf{i}}{2} x_2 =\\
&&
\partial_1 \partial_2 + \partial_1 \circ (\mathbf{i} x_1 +
z_2^2)^{1/2} + \frac{\mathbf{i}}{2} x_2 (\mathbf{i} x_1 +
z_2^2)^{-1/2}\partial_2 .
\een
Comparing the coefficients in front of $(\mathbf{i} x_1 +
z_2^2)^{-k/2-1/2}$ we get the following recursion relation 
\beq\label{rec-reln-2}
(-k+1) \partial_2 \psi^{(2)}_{k-2}(x_2) + (x_2 \partial_2
-k+1)\psi^{(2)}_{k-1}(x_2) =0. 
\eeq
Writing each $\psi^{(2)}_{k}(x_2) = \sum_{a=0}^\infty \psi^{(2)}_{k,a}
x_2^a$ and comparing the coefficients in front of $x_2^a$ we get 
\ben
(-k+1) (a+1) \psi^{(2)}_{k-2,a+1} + (a-k+1) \psi^{(2)}_{k-1,a} =0.
\een
Using these relations we get that $\psi^{(2)}_{k,a}=0$ for all $a>k+1$
and that 
$$
\psi^{(2)}_{k,a}=\frac{(a-k-2)\cdots (a-k-a-1)}{(k+1)\cdots (k+a) a!}\, 
\psi^{(2)}_{k+a,0}.
$$
Therefore the set of functions $\psi^{(2)}_{k}(x_2)$ $(k\geq 0)$ is
uniquely determined  from its restriction to $x_2=0$ as claimed. 

It remains to prove that the restriction $\Psi(x_1,0,z)$ is 
determined from the 1st equation. To be more precise, differentiating
the third equation in Theorem \ref{t2} with respect to $x_2$ and
restricting to $x_2=0$ 
we get 
\ben
\left. \Big( \partial_1\partial_2^2 \Psi\Big) \right|_{x_2=0} =
\frac{\mathbf{i}}{2}\, \Psi|_{x_2=0}.
\een 
Differentiating the first equation in Theorem \ref{t2} in $x_1$,
restricting to $x_2=0$, 
and using the above relation we get the following scalar differential 
equation for $\Psi|_{x_2=0}$
\ben
\Big(
\partial_1^{h+1} -(\mathbf{i} x_1 + z_1^h,\mathbf{i} x_1 + z_2^2)) \partial_1
-\frac{\mathbf{i}}{2} \Big) \Psi|_{x_2=0} =0.
\een
Conjugating with the exponential factors of the wave function, we get
that the series 
\ben
\sum_{k=0}^\infty \psi^{(1)}_k(0) (\mathbf{i} x_1 +
z_1^h)^{-k/h-1/2},\quad \psi^{(1)}_0(0):=1
\een
and 
\ben
\sum_{k=0}^\infty \psi^{(2)}_k(0) (\mathbf{i} x_1 +
z_2^2)^{-k/2-1/2},\quad \psi^{(2)}_0(0):=1,
\een
are anihilated by the differential operators
\beq\label{diff-op-1}
(\partial_1+ (\mathbf{i} x_1 + z_1^h)^{1/h})^{h+1} 
-(\mathbf{i} x_1 + z_1^h) (\partial_1+ (\mathbf{i} x_1 + z_1^h)^{1/h})
-\frac{\mathbf{i}}{2} 
\eeq
and 
\beq\label{diff-op-2}
\partial_1^{h+1} -(\mathbf{i} x_1 + z_2^2) \partial_1
-\frac{\mathbf{i}}{2},
\eeq
respectively. 
We claim that comparing the coefficients in front of $(\mathbf{i} x_1
+ z_1^h)^{-k/h-1/2}$ and $(\mathbf{i} x_1
+ z_2^2)^{-k/2-1/2}$ yields a recursion that uniquely determines
$\psi^{(a)}_k(0)$ for all $a=1,2$ and $k>0$.

Let us change the variables in the differential operator
\eqref{diff-op-1} via $y_1=(\mathbf{i} x_1 + z_1^h)^{1/h}$. Since 
\ben
\partial_1+ (\mathbf{i} x_1 + z_1^h)^{1/h} = y_1^{-h} ( y_1^{h+1} + \mathfrak{D}_1),
\een
where $\mathfrak{D}_1:=\frac{\mathbf{i}}{h} y_1 \partial_{y_1}$ we get that the
differential operator \eqref{diff-op-1} takes the form
\ben
y_1^{-h(h+1)} 
(y_1^{h+1}+\mathfrak{D}_1-\mathbf{i} \cdot h)\cdots
(y_1^{h+1}+\mathfrak{D}_1-\mathbf{i}\cdot 0) -\mathfrak{D}_1-y_1^{h+1} -\frac{\mathbf{i}}{2},
\een
where we used that $\mathfrak{D}_1 y_1^m= y_1^m(\mathfrak{D}_1+\mathbf{i} m/h)$. The above
operator can be expanded as a Laurent polynomial in $y_1^{h+1}$ with
coefficients (written on the left of the powers of $y_1$) that are
polynomials in $\mathfrak{D}_1$. Note that the coefficient in front of
$y_1^{h+1}$ is $0$, while the free term is
\ben
\sum_{k=0}^h (\mathfrak{D}_1-\mathbf{i} k) -
\mathfrak{D}_1- \frac{\mathbf{i}}{2}=
h\Big(\mathfrak{D}_1+\frac{\mathbf{i}}{2}\Big).
\een
Therefore the differential operator \eqref{diff-op-1} has the form
\ben
h\Big(\mathfrak{D}_1+\frac{\mathbf{i}}{2}\Big)+ \sum_{s=1}^h y_1^{-s(h+1)} \, P_s(\mathfrak{D}_1) ,
\een
where $P_s(\mathfrak{D}_1)\in \CC[\mathfrak{D}_1]$ are some polynomials. The above operator
anihilates the series $\sum_{k\geq 0} \psi^{(1)}_k(0)
y_1^{-k-N+1}$. Comparing the coefficients in front of $y_1^{-k-N+1}$,
we get the recursion relation 
\ben
-\mathbf{i} k \psi^{(1)}_k(0) + 
\sum_{s=1}^h \psi^{(1)}_{k-s(h+1)} (0) \,
P_s(-\mathbf{i}(k/h+1/2)) = 0
\een
valid for all $k\in \ZZ$ provided that we define
$\psi^{(1)}_\ell(0):=0$ for $\ell<0$. By definition
$\psi^{(1)}_0(0)=1$, so the above recursion 
determines uniquely $\psi^{(1)}_k(0)$ for all $k\geq 0$. 

The argument for the differential operator \eqref{diff-op-2} is
similar. After the change $y_2=(\mathbf{i} x_1 +z_2^2)^{1/2}$, the
differential operator \eqref{diff-op-2} takes the form
\ben
y_2^{-2(h+1)} (\mathfrak{D}_2-\mathbf{i} \cdot h)\cdots (\mathfrak{D}_2-\mathbf{i}\cdot 0) -
\mathfrak{D}_2-\frac{\mathbf{i}}{2},
\een
where $\mathfrak{D}_2=\frac{\mathbf{i}}{2} y_2 \partial_{y_2}$. 
This operator anihilates the series $\sum_{k\geq 0} \psi^{(2)}_k(0)
y_2^{-k-1}$. Comparing the coefficients in front of $y_2^{-k-1}$, we
get the following recursion
\ben
\frac{\mathbf{i} k}{2} \psi^{(2)}_k(0) +
(-1)^{N-1}\mathbf{i} \Big(\prod_{s=0}^h  (s-(k-1)/2)\Big) 
\psi^{(2)}_{k-2h-2}(0)=0.
\een
Note that for $0\leq k\leq 2h+1$ the second term is $0$, because
$\psi^{(2)}_{\ell}(0) $ is by definition $0$ for $\ell<0$.  Therefore
$\psi^{(2)}_{k}(0)=0$ for all $1\leq k\leq 2h+1.$ By definition
$\psi^{(2)}_{0}(0)=1$, so the above recursion uniquely determines
$\psi^{(2)}_{k}(0)$ for all $k>0$. This completes the proof of
Corollary \ref{c2}.

\section{Tau functions of Gaussian type}

The goal in this section is to prove Theorem \ref{t3}.
\subsection{Time evolution}

Suppose that $\Psi(x,z)$ is a wave function of the 2-component BKP hierarchy. We
would like to determine under what conditions the logarithm of the
corresponding tau-function is a quadratic form. 

By definition the components of the wave function have the form
\ben
\Psi^{(a)}(x,z_a) = \Big(1 +\sum_{k=1}^\infty w^{(a)}_k(x)
z_a^{-k}\Big) e^{x_a z_a}.
\een
The 2-component BKP hierarchy can be formulated as a deformation of the wave
function of the form
\ben
\Psi^{(a)} (\mathbf{t},z) =  \Big(1 +\sum_{k=1}^\infty w^{(a)}_k(\mathbf{t})
z_a^{-k}\Big) e^{\sum_{k\in \ZZ_{\rm odd}^+} t^a_k z_a^k},
\een
where as usual we identify $t^1_1=x_1$ and $t^2_1=x_2$. 
The pseudodifferential operators
\ben
S^{(a)}(\mathbf{t},\partial_a) = 1 +\sum_{k=1}^\infty w^{(a)}_k(\mathbf{t})
\partial_a^{-k},\quad \partial_a=\partial/\partial x_a\quad (a=1,2)
\een
and
\ben
\L_a(\mathbf{t},\partial_a):=
S^{(a)}(\mathbf{t},\partial_a) \partial_a
S^{(a)}(\mathbf{t},\partial_a) ^{-1} \quad (a=1,2)
\een
are called {\em dressing operators} and {\em Lax
  operators}, respectively.   The dressing operators are called {\em
  wave operators} if the following system of differential equations is
satisfied
\ben
\partial_{t^a_k} \Psi(\mathbf{t},z) =
(\L_a(\mathbf{t},\partial_a) ^k)_+ \, \Psi(\mathbf{t},z),\quad
a=1,2,\quad k\in \ZZ_{\rm odd}^+,
\een
where the $+$ index of a pseudodifferential operator means truncating
the non-differential operator part. 
The above system of differential equations uniquely determines the wave
operators from their restriction to $t^1_1=x_1,t^2_1=x_2$, and
$t^a_m=0$ for $m>1$. In particular, the deformation
$\Psi(\mathbf{t},z)$ of the wave 
function is uniquely determined and it is still refered to as a wave
function. It is a very non-trivial result that for every wave function
$\Psi(\mathbf{t},z)$ there exists a formal series $\tau(\mathbf{t})$,
called {\em tau function}, such that 
\ben
\Psi^{(a)}(\mathbf{t},z_a) =
\frac{\Gamma(\mathbf{t}^a, z_a)\tau(\mathbf{t})}{\tau(\mathbf{t})}.
\een
Such a formal series is unique up to a constant factor and it
satisfies the Hirota bilinear equations of the 2-component BKP hierarchy. 
It is straightforward to prove that $\log \tau$ has at most quadratic
terms in $\mathbf{t}$ if an only if the non-deformed Lax operators
$\L_a(x,\partial_a)= \L_a(0,\partial_a)$ are constant, i.e.,
independent of $x$. 

\subsection{Properties}
Suppose that $\tau(\mathbf{t})$ is a tau-function of
the Kac--Wakimoto hierarchy that has the form
\eqref{tau-gaussian}. Let us assume also that the coefficients
$W^{ab}_{k\ell}$ are symmetric: $W^{ab}_{k\ell}=W^{ba}_{\ell k}$. 
The components of the corresponding wave function have the form
\begin{align}
\notag
\Psi^{(1)}(\mathbf{t},z_1) = &  \Psi^{(1)}(0,z_1) \times\\
\notag
&
\exp\Big(
\sum_{k\in \ZZ^+_{\rm odd} } 
t^1_k \Big(z_1^k-2\sum_{\ell\in \ZZ^+_{\rm odd} }  W^{11}_{k\ell} \frac{z_1^{-\ell}}{\ell}\Big) +
t^2_k \Big(-2\sum_{\ell\in \ZZ^+_{\rm odd} }  W^{21}_{k\ell} \frac{z_1^{-\ell}}{\ell}\Big) \Big)
\end{align}
and 
\begin{align}
\notag
\Psi^{(2)}(\mathbf{t},z_2) = & \Psi^{(2)}(0,z_2) \times \\
\notag
&
\exp\Big(
\sum_{k\in \ZZ^+_{\rm odd} } 
t^1_k \Big(
-2\sum_{\ell\in \ZZ^+_{\rm odd} }  W^{12}_{k\ell} \frac{z_2^{-\ell}}{\ell}\Big) +
t^2_k \Big(
z_2^k-2\sum_{\ell\in \ZZ^+_{\rm odd} }  W^{22}_{k\ell}
  \frac{z_2^{-\ell}}{\ell}\Big)
\Big) ,
\end{align}
where 
\ben
 \Psi^{(a)}(0,z_a) = \exp\Big(
2\sum_{k,\ell\in  \ZZ_{\rm odd}^+} W^{aa}_{k\ell}
\frac{z_a^{-k-\ell}}{k\ell}
\Big),\quad a=1,2.
\een
The first observation is that the subspace $U\in
\operatorname{Gr}_2^{I,(0)}$ corresponding to the wave function is a
$\CC[X_1,X_2]$-module, where $X_a=X_a(z)\in V$ ($a=1,2)$ are defined by  
$\partial_{t^a_1} \Psi = X_a(z) \Psi$, i.e., 
\beq\label{x1}
X_1=\Big(z_1-2\sum_{\ell\in \ZZ_{\rm odd}^+} W^{11}_{1\ell}
z_1^{-\ell}/\ell\ ,\ 
-2\sum_{\ell\in \ZZ_{\rm odd}^+} W^{12}_{1\ell}
z_2^{-\ell}/\ell
\Big)
\eeq
and 
\beq\label{x2}
X_2=\Big(-2\sum_{\ell\in \ZZ_{\rm odd}^+} W^{21}_{1\ell}
z_1^{-\ell}/\ell\ ,\ 
z_2-2\sum_{\ell\in \ZZ_{\rm odd}^+} W^{22}_{1\ell}
z_2^{-\ell}/\ell
\Big),
\eeq
where the action of $X_a$ on the wave function is via componentwise
multiplication. Note also that $U=\CC[X_1,X_2] \Psi(0,z)$.
\begin{lemma}\label{le:gaussian-1}
a) Put 
$$
t_N:=2\partial_1\partial_2 \log \tau=
2W^{12}_{11}=2W^{21}_{11},
$$ 
then $X_1X_2+t_N=0$.

b) There are parameters $t_i$ ($1\leq i\leq N-1$), such that 
\ben
X_1^h + \sum_{i=1}^{N-1} t_i X_1^{h-2i} + X_2^2 = (z_1^h,z_2^2).
\een 

c) There are uniquely defined polynomials 
$$
p_k^{(a)} (t,x)\in \CC[t,x],\quad k\in \ZZ_{\rm odd}^+,\quad a=1,2, 
$$
where $t=(t_1,t_2,\dots,t_N)$ such that with respect to $x$ they are
monic of degree $k$ and 
\ben
p_k^{(a)}(t,X_a(z)) = z_a^k e_a + O(z_1^{-1}) e_1 + O(z_2^{-1}) e_2.
\een
\end{lemma}
\proof
a) Recall property $(W2)$ of the wave function (see
Section \ref{sec:2-BKP}). There exists function $q(x_1,x_2)$ such
that $(\partial_1\partial_2+q)\Psi=0$. Comparing the coefficients in
front of $z_1^0$ we get that $q=2\partial_1\partial_2 \log
\tau=t_N$. By definition $\partial_1\partial_2 \Psi = X_1 X_2\Psi$, so
the relation $X_1X_2+t_N=0$ follows.

b) Note that 
\ben
X_1^{2m} = z_1^{2m}e_1 + (O(z_1^{2m-2}), O(z_2^{-2m}))
\een
and 
\ben
X_2^2 = z_2^2 e_2 + (O(z_1^{-2}),O(z_2^{0})).
\een
Therefore we can choose the coefficients $t_1,\dots, t_{N-1}$ in such
a way that 
\ben
F:=X_1^h + \sum_{i=1}^{N-1} t_i X_1^{h-2i} + X_2^2 -
(z_1^h,z_2^2)
\een
satisfies $F (e_1+\mathbf{i} e_2)\in V_0$. 
This however implies that $F \Psi(0,z)\in U\cap V_0=\{0\}$, i.e.,
$F=0$.

c) Using part a) and b) we can express all $W^{ab}_{1\ell}$,
($a,b=1,2$, $\ell \in \ZZ^+_{\rm odd}$) as polynomials in $t$. It
follows that if $k\in  \ZZ^+_{\rm odd}$, then 
\ben
X_a(z)^k = z_a^k e_a + O(z_a^{k-2}) e_a +  O(z_2^{-k}) e_b, 
\een
where $b$ is such that $\{a,b\}=\{1,2\}$, the RHS is a Laurent series
in $z_1$ and $z_2$ that involves only 
odd powers of $z_1$ and $z_2$ and its coefficients are polynomials in
$t$. Clearly we can remove the lower order terms on the RHS  that
involve non-negative powers of $z_a$ by adding a linear combination of
$X_a(z)^{k-2},\dots, X_a(z)$.
\qed

We claim that the coefficients $W^{ab}_{k\ell}$ are polynomials in the
parameters $t$. We already argued that part a) and b) of Lemma
\ref{le:gaussian-1} implies that $W^{ab}_{1\ell}$ are polynomials in
$t$. Let us prove that $W^{1b}_{k\ell}$ is polynomial in $t$. The
argument for $W^{2b}_{k\ell}$ is similar. 
Put 
$$
A=\Big(z_1^k-2\sum_{\ell\in \ZZ^+_{\rm odd} }  W^{11}_{k\ell}
\frac{z_1^{-\ell}}{\ell}\Big)e_1 + 
\Big(
-2\sum_{\ell\in \ZZ^+_{\rm odd} }  W^{12}_{k\ell}
\frac{z_2^{-\ell}}{\ell}\Big) e_2.
$$
By definition 
$
\partial_{t^1_k} \Psi= A\Psi,
$
so $A\Psi(0,z)\in U$. On the other hand
$$
(A-p_k^{(1)}(t,X_1))\Psi(0,z)\in U\cap V_0=\{0\}.
$$ 
Hence $A=p_k^{(1)}(t,X_1)$. 

\begin{lemma}\label{le:gaussian-2}
The components of the wave function $\Psi(0,z)$ are given by the
following formulas
\ben
\Psi^{(a)}(0,z_a)^2 = \frac{z_a\partial_{z_a} X_{aa}}{X_{aa}},\quad a=1,2,
\een
where we define $X_{ab}$ by  $X_a= X_{a1}e_1+X_{a2} e_2$. 
\end{lemma}
\proof
Since $U$ is an isotropic subspace the vectors $X_1^m \Psi(0,z)$ must
be isotropic for all integers $m\geq 0$. This imples that 
\ben
\operatorname{Res}_{z_1=0}
X_{11}^{2m}\Psi^{(1)}(0,z_1)^2\frac{dz_1}{z_1} = \delta_{m,0} .
\een
Note that series $\Psi^{(1)}(0,z_1)^2=1+E_1 z_1^{-2} + E_2
z_1^{-4}+\cdots$ is uniquely determined in terms of $X_1$ from the above residue
relations: the relation for $m=1$ determines $E_1$, the relation for $m=2$ determines
$E_2$, etc.. On the other hand
\ben
\operatorname{Res}_{z_1=0} X_{11}^{2m}\frac{dX_{11} }{X_{11}} = 
-\operatorname{Res}_{X_{11}=\infty} X_{11}^{2m}\frac{dX_{11} }{X_{11} } = \delta_{m,0}.
\een
This relation implies the formula for $a=1$. The proof of the formula
for $a=2$ is similar.
\qed

\subsection{Proof of Theorem \ref{t3}}

We have constructed a map from the set of tau-functions of the
Kac--Wakimoto hierarchy of Gaussian type to $\CC^N$. The map is given by the
parameters $t=(t_1,\dots,t_N)$ from Lemma
\ref{le:gaussian-1}. Moreover, we proved that the coefficients
$W^{ab}_{k\ell}$ can be expressed as polynomials in $t$. This
proves that the map is injective. It remains only to prove that it is
surjective. 

Given parameters $t=(t_1,\dots,t_N)$ we can uniquely
define $W^{ab}_{1\ell}$ such that $X_1$ and $X_2$ given by formulas
\eqref{x1} and \eqref{x2} respectively are solutions to the polynomial
equations in parts a) and b) of Lemma \ref{le:gaussian-1}. Furthermore, we
define the remaining coefficients $W^{ab}_{k\ell}$ via the relations
\ben
\Big(z_1^k-2\sum_{\ell\in \ZZ^+_{\rm odd} }  W^{11}_{k\ell}
\frac{z_1^{-\ell}}{\ell}\Big)e_1 + 
\Big(
-2\sum_{\ell\in \ZZ^+_{\rm odd} }  W^{12}_{k\ell}
\frac{z_2^{-\ell}}{\ell}\Big) e_2:= p_k^{(1)}(t,X_1)
\een 
and 
\ben
\Big(-2\sum_{\ell\in \ZZ^+_{\rm odd} }  W^{21}_{k\ell}
\frac{z_1^{-\ell}}{\ell}\Big)e_1 + 
\Big(
z_2^k-2\sum_{\ell\in \ZZ^+_{\rm odd} }  W^{22}_{k\ell}
\frac{z_2^{-\ell}}{\ell}\Big) e_2:= p_k^{(2)}(t,X_2),
\een 
where $p_k^{(a)}(t,x)$ are the polynomials from Lemma
\ref{le:gaussian-1}, part c).  We claim that 
\ben
\Psi^{(a)}(\mathbf{t},z_a)=
\Big( \frac{z_a\partial_{z_a} X_{aa}}{X_{aa}}\Big)^{1/2} \ 
\exp\Big(
\sum_{k\in \ZZ_{\rm odd}^+} 
t^1_k p_k^{(1)}(t,X_{1a}) +
t^2_k p_k^{(2)}(t,X_{2a})
\Big)  
\een
are the components of a wave function. 
By definition
\ben
\partial_a \Psi^{(a)} = X_{aa}(z_a) \Psi^{(a)}=X_{aa}(\L_a) \Psi^{(a)}.
\een
Hence the Lax operators can be found by solving the equations
\ben
\partial_a = \L_a -2\sum_{\ell\in \ZZ_{\rm odd}^+}
W^{aa}_{1\ell}\frac{\L_a^{-\ell}}{\ell} \quad (a=1,2)
\een
for $\L_a$ in terms of $\partial_a$. In particular, the Lax operators
are independent of all dynamical variables $\mathbf{t}$. Note that by
definition 
\ben
p_k^{(a)}(t,X_{aa}(\L_a)) = \L_a^k + O(\L_a^{-1}).
\een
Taking the differential operator part of both sides and recalling
$X_{aa}(\L_a)=\partial_a$ we get that 
\ben
(\L_a(0,\partial_a)^k)_+ = p_k^{(a)}(t,\partial_a).
\een
On the other hand, by definition 
$$
\partial_a \Psi^{(b)}(\mathbf{t},z_b)=X_{ab}(z_b)
\Psi^{(b)}(\mathbf{t},z_b).
$$
Therefore
\begin{align}
\notag
\partial_{t^a_k}\Psi^{(b)}(\mathbf{t},z_b) & =
p_k^{(a)}(t,X_{ab}) \Psi^{(b)} (\mathbf{t},z_b) =
p_k^{(a)}(t,\partial_a) \Psi ^{(b)}(\mathbf{t},z_b) =\\
\notag
& =
(\L_a(0,\partial_a)^k)_+\Psi^{(b)}(\mathbf{t},z_b) .
\end{align}
This proves that $\Psi(\mathbf{t},z)$ satisfies the differential
equations of the 2-component BKP hierarchy. We need only to varify that
$\Psi(x,z)$ is a wave-function, i.e., it satisfies the axioms
(W1)--(W3). Axioms (W1) and (W2) hold by definition. Only (W3)
requires an argument. Recalling the definition of $\Psi(x,z)$, we get
that the pairing
\ben
(\Psi(x',z), \Psi(x'',z))
\een
is a difference of two residues 
\beq\label{res-z1}
\operatorname{Res}_{z_1=0} \frac{d X_{11}}{X_{11}}
e^{(x_1'-x_1'')X_{11} + (x_2'-x_2'') X_{21}} 
\eeq
and 
\beq\label{res-z2}
\operatorname{Res}_{z_2=0} \frac{d X_{22}}{X_{22}}
e^{(x_1'-x_1'')X_{12} + (x_2'-x_2'') X_{22}} .
\eeq
Note that $X_1X_2=-t_N$ implies that $X_{11}X_{21}=-t_N$ and
$X_{12}X_{22}=-t_N$. Since $X_{11}$ is a coordinate near $z_1=\infty$,
the first residue can be written as 
\ben
-\operatorname{Res}_{X_{11}=\infty} \frac{d X_{11}}{X_{11}}
e^{(x_1'-x_1'')X_{11} + (x_2'-x_2'') (-t_N/X_{11})}, 
\een
where the 1-form should be expanded as a power series in $x_1'-x_1''$
and $x_2'-x_2''$. The residue is interpreted first formally as the
negative of the coefficient in front of $dX_{11}/X_{11}$. Note
however that the coefficients of the formal power series are Laurent
polynomials in $X_{11}$, so each residue can be interpreted
analytically. Similarly the second residue \eqref{res-z2} equals
\ben
-\operatorname{Res}_{X_{22}=\infty} \frac{d X_{22}}{X_{22}}
e^{(x_1'-x_1'')(-t_N/X_{22}) + (x_2'-x_2'') X_{22}}.
\een
Changing the variables $X_{22}\mapsto -t_N/X_{11}$ yields
\ben
\operatorname{Res}_{X_{11}=0} \frac{d X_{11}}{X_{11}}
e^{(x_1'-x_1'')X_{11} + (x_2'-x_2'') (-t_N/X_{11})},
\een
so we need to prove that 
\ben
(\operatorname{Res}_{X_{11}=0} + \operatorname{Res}_{X_{11}=\infty} )\frac{d X_{11}}{X_{11}}
e^{(x_1'-x_1'')X_{11} + (x_2'-x_2'') (-t_N/X_{11})} =0.
\een
This however follows from the residue theorem for $\mathbb{P}^1$, because the residues
involved are residues of Laurent polynomials in $X_{11}$, so the only
poles are at $X_{11}=0$ and $X_{11}=\infty.$ This completes the proof
that $\Psi(\mathbf{t},z)$ is a wave function of the 2-component
BKP hierarchy. The corresponding point $U\in \operatorname{Gr}_2^{I,(0)}$ is
invariant under multiplication by $X_1$ and $X_2$ and since by
definition $X_1$ and $X_2$ satisfy the relation in part b) of Lemma
\ref{le:gaussian-1} we get that $U$ is invariant under multiploication
by $(z_1^h,z_2^2)$. Hence $\Psi(\mathbf{t},z)$ is a wave function of
the Kac--Wakimoto hierarchy. Furthermore, since the Lax operators are
constant the tau-function corresponding to $\Psi(\mathbf{t},z)$ must
have the form 
\ben
\widetilde{\tau}(\mathbf{t}) = \exp\Big(
\sum_{a=1,2}\sum_{k\in \ZZ_{\rm odd}^+}
\widetilde{W}^a_k t^a_k + \frac{1}{2}
\sum_{a,b=1,2}\sum_{k,\ell\in \ZZ_{\rm odd}^+}
\widetilde{W}^{ab}_{k\ell} t^a_k t^b_\ell
\Big). 
\een
We need only to prove that $\widetilde{W}_k^a=0$ and
$\widetilde{W}^{ab}_{k\ell}=W^{ab}_{k\ell}$. Writing
$\Psi(\mathbf{t},z)$ in terms of the tau-function 
$\widetilde{\tau}$ we get that $\log \Psi^{(a)}(\mathbf{t},z_a)$ can be
written as a power series in $\mathbf{t}$ that involves at most linear
terms. Recalling the definition of $\Psi(\mathbf{t},z)$ and comparing
the free terms, we get  
\ben
-2\sum_{k\in \ZZ_{\rm odd}^+} \widetilde{W}_k^a \frac{z_a^{-k}}{k} 
+2\sum_{k,\ell \in \ZZ_{\rm odd}^+} \widetilde{W}_{k\ell}^{aa}
\frac{z_a^{-k-\ell}}{k\ell} =
\frac{1}{2} \log  \frac{z_a\partial_{z_a} X_{aa}}{X_{aa}},\quad a=1,2.
\een
The RHS involves only even powers of $z_a^{-1}$. Therefore all
coefficients $\widetilde{W}_k^a =0$. Comparing the coefficients in
front of $t^a_k$ in $\log \Psi^{(b)}(\mathbf{t},z_b)$ we get
\ben
\sum_{\ell\in \ZZ_{\rm odd}^+}
\widetilde{W}^{ab}_{k\ell}\frac{z_b^{-\ell}}{\ell} = 
\sum_{\ell\in \ZZ_{\rm odd}^+}
W^{ab}_{k\ell}\frac{z_b^{-\ell}}{\ell}  ,\quad a,b=1,2,\quad k\in
\ZZ_{\rm odd}^+.
\een
This implies that $\widetilde{W}^{ab}_{k\ell}= W^{ab}_{k\ell}$.
\qed

\bibliographystyle{amsalpha}

\end{document}